\documentclass{article}
\usepackage[utf8]{inputenc} 
\usepackage{amsmath, amssymb, amsthm} 
\usepackage{graphicx} 
\usepackage{hyperref} 
\usepackage{geometry} 
\usepackage{setspace} 
\usepackage{booktabs}
\usepackage{enumitem}
\usepackage{setspace}

\usepackage[table]{xcolor}
\usepackage[utf8]{inputenc}
\usepackage[english]{babel}
\usepackage{amsmath}
\usepackage{graphicx}
\usepackage{amssymb}
\usepackage[table]{xcolor}

\usepackage{caption}
\usepackage{subcaption}
\usepackage{textcomp}
\usepackage{algorithm}
\usepackage{algpseudocode}
\usepackage{graphicx}%
\usepackage{multirow}%
\usepackage{amsmath,amssymb,amsfonts}%
\usepackage{amsthm}%
\usepackage{mathrsfs}%
\usepackage[title]{appendix}%
\usepackage{xcolor}%
\usepackage{textcomp}%
\usepackage{manyfoot}%
\usepackage{booktabs}%
\usepackage{algorithm}%
\usepackage{algorithmicx}%
\usepackage{algpseudocode}%
\usepackage{listings}%
\usepackage{geometry}
\title{Computing the steady-state probabilities of a tandem queueing system, an ML approach}
\author{Eliran Sherzer \\
    Ariel University, Israel \\
    \texttt{eliransh@ariel.ac.il} }

\newtheorem{remark}{Remark}%

\begin{document}

\maketitle

\begin{abstract}
Tandem queueing networks are widely used to model systems where services are provided in sequential stages. In this study, we assume that each station in the tandem system operates under a general renewal process. Additionally, we assume that the arrival process for the first station is governed by a general renewal process, which implies that arrivals at subsequent stations will likely deviate from a renewal pattern.

This study leverages neural networks (NNs) to approximate the steady-state distribution of the marginal number of customers at each station in the tandem queueing system, based on the external inter-arrival and service time distributions.

Our approach involves decomposing each station and estimating the departure process by characterizing its first five moments and auto-correlation values, without limiting the analysis to linear or first-lag auto-correlation. We demonstrate that this method outperforms existing models, establishing it as state-of-the-art.

Furthermore, we present a detailed analysis of the impact of the $i^{th}$ moments of inter-arrival and service times on steady-state probabilities, showing that the first five moments are nearly sufficient to determine these probabilities. Similarly, we analyze the influence of inter-arrival auto-correlation, revealing that the first two lags of the first- and second-degree polynomial auto-correlation values almost completely determine the steady-state probabilities of a $G/GI/1$ queue.
\end{abstract}

\section{Introduction} \label{sec:intro}

Tandem queueing networks are extensively utilized to model systems where services are provided in sequential stages. In such systems, generic customers enter a multi-stage service facility and proceed through each service stage in succession. This modeling approach is applicable across various fields, including transportation, manufacturing, logistics, computer networking, telecommunications, and numerous everyday service operations. Examples include production lines with different processing steps handled by various machines, airplane maintenance and refueling, the ordering and delivery process of goods or services, and routing messages or data packets through multiple hops in the Internet or other communication networks.

Service operation policies frequently depend on quantitative descriptions of system performance, known as performance measures, such as waiting time, queue length, and workload within the system. Effective decision-making for service operations requires an accurate understanding of these performance measures.

A common approach to analyzing the performance of complex queueing models is through computer simulation (e.g., see \cite{doi:10.1080/07408170590899625, 10.1145/2000494.2000497}. However, as highlighted by \cite{doi:10.1287/opre.2016.1554}, a major drawback of simulation-based optimization methods is the often excessive computation time needed to achieve optimal solutions for service operation problems involving multidimensional stochastic networks.

The class of queueing networks that can be solved analytically requires stringent assumptions that are rarely met, while more realistic models are exceedingly difficult to analyze precisely. Alternatively,  one may use approximations analytic methods; however, such models are often unreliable. Therefore, a different approach is needed. This paper presents a steady-state analysis of a tandem queueing system using machine learning (ML).

\begin{figure}
\centering
\includegraphics[scale=0.42]{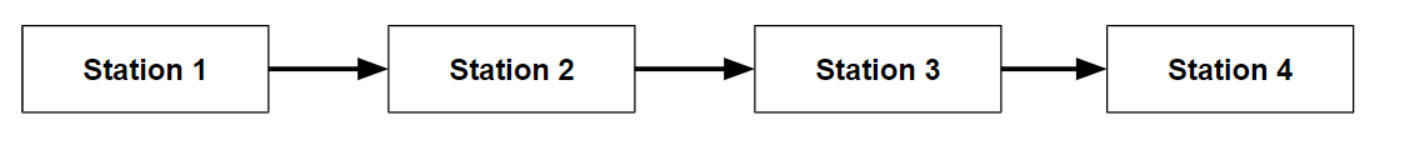}
\caption{An architecture of a tandem queueing system with four server stations. }
\label{fig:tandem_archi}
\end{figure}

This study assumes that services in each station in the tandem system follow a general renewal process. It is assumed that the arrival process for the first station (i.e., the external arrival) is also a general renewal process. Of course, as a result, the arrivals for the following stations are expected to be non-renewal. For example, in Figure~\ref{fig:tandem_archi}, the arrival processes for stations 2, 3, and 4 are non-renewal since the departure process from stations 1, 2, and 3 are non-renewal, respectively.
This study uses neural networks (NNs) to approximate the marginal number of customers' steady-state distribution in each station of an $m \geq 2$ tandem queueing system, given the external inter-arrival distribution and the service time distribution.   
We note that we do not attempt to approximate the steady-state distribution of the first station (i.e., station 1 in Figure~\ref{fig:tandem_archi}), that is of a $GI/GI/1$ system, which was done in a previous study~\cite{sherzer23}. There, the steady-state distribution was approximated via an NN given the first five moments of the inter-arrival and service distribution, with an outstanding level of accuracy.

As typically done in the literature when analyzing queueing systems, we decompose each station and analyze it separately \cite{Whitt1994, 6770716, doi:10.1287/mnsc.34.1.75}. Generally speaking, the decomposition approach under a tandem system operates as follows: for each station, we evaluate the departure process and use it as the arrival process for the next station. Given both the arrival (possibly non-renewal) and the service times, we evaluate the steady-state number of customers in the system. 

We approximate the steady-state number of customers at each station via the arrival process and service times. Since departures do not need to follow a renewal process, it is non-trivial which information should be fed to the NN that approximates the steady-state probabilities of the following station. 
That is, the departure process, and hence the arrival process of the next station, consists of two parts. The first part is the equilibrium departure distribution, and the second is the auto-correlation function, which reflects dependencies between departures.   

In previous studies, such information was given as the departure via an approximation of the squared coefficient of variation (SCV) (e.g.,~\cite{Whitt1994, 6770716, SAGRON2017161}), along with the first moment, which is known. However, the SCV and the first moment do not reflect the dependencies between arrivals, which, according to~\cite{CIVELEK20211031}, is essential for such computations. As such, an alternative approach was given by~\cite{https://doi.org/10.1002/nav.22010}, which the defined index of dispersion for counts (IDC)~\cite{WHITT201999}: "The IDC is a scaled version of the variance–time function. In particular, if $A(t)$ is an arrival counting process assumed to be stationary with rate $\lambda$,  the IDC is:

\begin{align*}
    I_c(t)\equiv \frac{Var(A(t))}{E[A(t)]}= \frac{Var(A(t))}{\lambda t}, \thinspace t \geq 0,
\end{align*}
\noindent where $\equiv$ denotes equality by definition."

The IDC was proven very efficient and outperformed other methods based on SCV approximation when comparing the steady-state analysis of the next station~\cite{https://doi.org/10.1002/nav.22010}, most likely because the IDC holds information regarding the dependencies between arrivals.   Yet, despite the improvement, significant errors can still be made, in many cases larger than 10\%~\cite{WHITT201999}.

In contrast, we propose a different way to present the departure process. We propose combining the inter-departure first $n\in \mathbb{N}$ moments with their corresponding auto-correlation values. By auto-correlation values, we do not refer to only the $1^{st}$-lag linear auto-correlation. Instead, we feed the network values concerning the first $n_1 \in \mathbb{N}$-lag auto-correlation. For each lag $k \leq n_1$, we feed polynomial auto-correlation values up to a predefined degree $n_2 \in \mathbb{N}$.  In Section~\ref{sec:prob_formulation}, we give an exact definition for all auto-correlation values.

$n$, $n_1$, and $n_2$ are parameters which determine the input of our model. We wish to feed the network input to produce the most accurate prediction. Values of $n$, $n_1$, and $n_2$ which are too small will not hold enough information; for example, if we feed only the first inter-arrival service time moment (i.e., $n=1$), ignoring the second, third, and so forth, will remove vital information to capture the dynamic nature of the queue according to~\cite{CIVELEK20211031}. In contrast,  $n$, $n_1$, and $n_2$ that are too large (i.e., feeding the network with very high moments, high-lags correlation values with higher polynomial levels) may not include vital information. Still, they simultaneously may insert undesired noise, which can cause an accuracy decrease~\cite{sherzer23}. 

As such, we conduct an empirical study that examines the model accuracy as a function of $n$, $n_1$, and $n_2$. Besides optimizing the model performances, the results of such an experiment can shed light on a broad queueing question. How do the moments of the renewal service time, the moments of the non-renewal arrival process, and its auto-correlation specification affect stationary queue length behavior? 

From~\cite{sherzer23}, which, as mentioned above, studied a  $GI/GI/1$ system, suggests that the effect of the $6^{th}$ moments of both the inter-arrival and service times is insignificant. Therefore, there was no gain in accuracy when using them as inputs. Their results also quantified the accuracy improvement for each moment.  Our experiment suggests that such results reserve from the transition to a $G/GI/1$ queue; hence, we use $n=5$. 

As for the value of $n_1$, We know from~\cite{CIVELEK20211031} that both first and second-lag linear auto-correlation values affect the queue dynamic of the following queue, suggesting that $n_1$ should be at least 2. Our empirical results indicate that $n_1$ is precisely 2. We did not find a similar study anywhere in the literature that discusses the effect of the inter-arrival non-linear auto-correlation on queue dynamics. Nevertheless, empirical results suggest that $n_2$ is also 2.

So far, we have focused on predicting steady-state probabilities given the non-renewal departure process moments and correlation values. However, obtaining them is also non-trivial. To our knowledge, no results in the literature accurately predict the departure moments and auto-correlation from both $GI/GI/1$ and $G/GI/1$ systems. The most relevant results are by \cite{math12091362}, which can accurately compute the departure moments of the $PH/GI/1$ system and its first and second lag linear auto-correlation. However, their results do not apply to a general $GI/GI/1$, and non-linear correlations cannot be derived. Furthermore, their results do not apply to a $G/GI/1$ system. Thus, to complete the analysis of the tandem queueing system, our approach also requires approximating the departure process of both  $GI/GI/1$ and $G/GI/1$ queues. Thus, we dedicate an NN to learning the departure process from $GI/GI/1$ and $G/GI/1$ queues, and their output will then be used to approximate the steady-state probabilities.  Next, we describe the full framework of our method.

\begin{figure}
\centering
\includegraphics[scale=0.42]{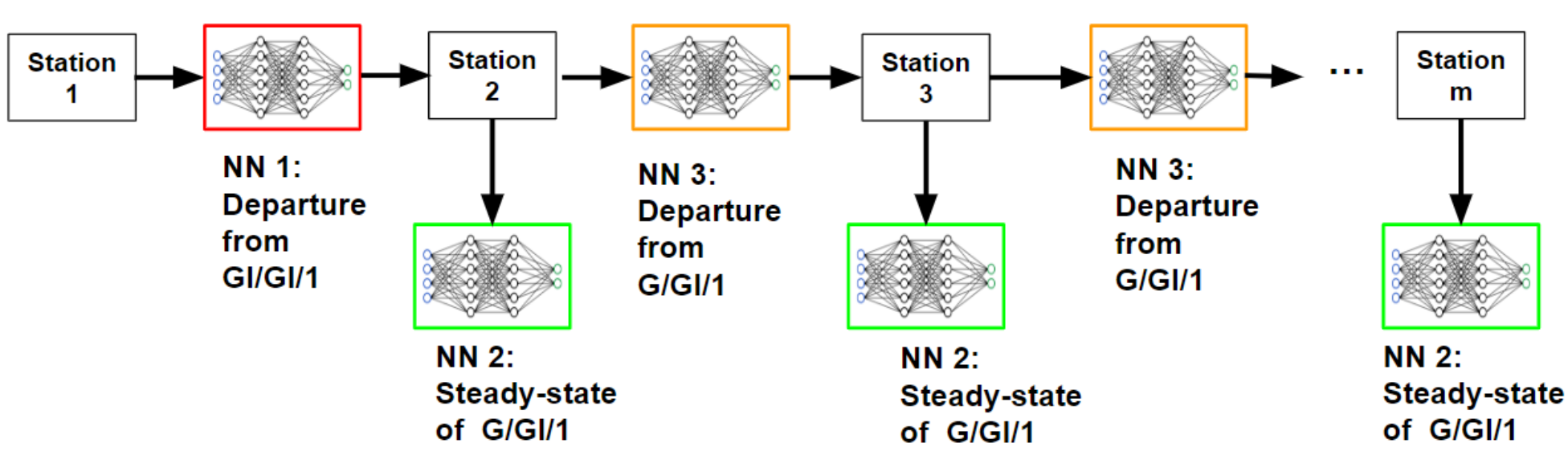}
\caption{Illustration of the solution overview of an $m$ station tandem system. }
\label{fig:tandem_archi_nn}
\end{figure}

\textbf{Solution overview:} To approximate the number of customers' steady-state distribution for each station $2 \leq j \leq m$, our method requires three different approximations, each of which has a designated neural network as depicted in Figure~\ref{fig:tandem_archi_nn}. The first is the departure moments and auto-correlation from the first station (namely station 1), which is a $GI/GI/1$ station. We denote this neural network by \textbf{NN 1} (marked by a red rectangle in Figure~\ref{fig:tandem_archi_nn}). The input of \textbf{NN 1} is the first $n$ moments of the external inter-arrival and service time distributions. The output is the departure process's first $n$ moments and auto-correlation values, consisting of the first $n_1$ lags and the first $n_2$ polynomial order auto-correlation values for each lag (which are defined in Section~\ref{sec:prob_formulation}).

Once the departure of station 1 is approximated, we estimate steady-state probabilities at the next station (namely station 2) with an NN we denote by  \textbf{NN 2} (marked by a green rectangle in Figure~\ref{fig:tandem_archi_nn}). The input is the service time first $n$ moments and the output of \textbf{NN 1} (i.e., the departure moments and auto-correlation value from station 1). The output is the steady-state probabilities of the number of customers at station 2. The next step is estimating the departure process of station 2, which is a $G/GI/1$ queue. This is done via  \textbf{NN 3} (marked by an orange rectangle in Figure~\ref{fig:tandem_archi_nn}). The input of \textbf{NN 3} is the departure of station 2, which is the output of  \textbf{NN 2} and the first $n$ service time moments of station 2. 

The next step estimates station 3 steady-state probabilities via \textbf{NN 2}. The input is the service time of station 3 and the departure process from station 2, estimated by \textbf{NN 3}—the output, similar to station 2, steady-state probabilities. We now hit a loop where we compute the departure process via \textbf{NN 3}, to be used as an input along with the service time moments at the next station, and compute the steady-state probabilities using \textbf{NN 2} until we reach station $m$. Notice that the number of NNs required is always three regardless of the size of $m$.

\begin{remark}
Note that \textbf{NN 1} is a special case of \textbf{NN 3}. We could use only \textbf{NN 3}, while for non-renewal arrival processes, we simply insert zeros for all auto-correlation entries. However, we believe the more natural way to build this framework is to increase the complexity gradually, as from the queueing perspective, \textbf{NN 1} and \textbf{NN 3} are fundamentally different. 
\end{remark}

To execute the framework described above and depicted in Figure~\ref{fig:tandem_archi_nn}, we need to train \textbf{NN 1}, \textbf{NN 2}, and \textbf{NN 3}. Since this is supervised learning, many instances of input and output are required for each NN. The challenge here is two-fold. The first is to generate many instances for each one of the queues for \textbf{NN 1}, \textbf{NN 2}, and \textbf{NN 3}, which are determined by the external inter-arrival and service time distributions. The second is labeling the input. Clearly, the training set must be diverse with accurate labeling to have a tool that accurately predicts a tandem stationary system.

To address the challenge of generating diverse training instances that capture a broad range of possible inputs, we build upon the work of~\cite{sherzer23}, who utilized the Phase-Type (PH) family of distributions. This family is known for its density with respect to all non-negative distributions, meaning that any given non-negative distribution can be approximated by a PH distribution to a specified level of accuracy (see Theorem 4.2, Chapter III~\cite{Asmussen2003}). In their work, \cite{sherzer23} developed an algorithm to generate varied PH instances, allowing their $GI/GI/1$ model to accurately predict steady-state occupancy distributions across a wide range of datasets. Thus, by using their PH sampler, we generate diverse sets of external inter-arrival and $m$ different service time distributions.

For the second challenge mentioned—labeling the model input—our only viable option is to use discrete event simulations. The primary drawback of simulation is its lengthy run time. At first glance, this may seem counterproductive: we turn to machine learning because simulations are time-consuming, yet to implement the machine learning approach, we need to perform numerous simulations to generate a training dataset.

However, the key advantage here is that labeling occurs offline, where it’s possible to temporarily allocate substantial computing resources, and response time is less critical. The machine learning approach is ideally suited for applications where fast response times are essential. In such cases, our approach allows for near real-time inference. In contrast, simulations can take anywhere from a few minutes to several hours. Additionally, our approach only requires a finite set of training examples. Once the model is trained, it can make inferences for infinite examples.

Furthermore, for labeling our training examples, we only need to simulate a 2-station tandem system, which allows us to do inference for any $m\geq 2$ size system.  The reason a 2-station is sufficient is that the departure from the first station can be used to train \textbf{NN 1}, the steady-state probabilities of the second station can be used for \textbf{NN 2}, and departures from the second station can be used for \textbf{NN 3}.  

We conduct a thorough performance evaluation to ensure our approach is accurate for the trained NNs. We show that our model is accurate on a diverse test set.  Furthermore, we compare our model to the existing approximation methods. We note that no model estimates the steady-state probabilities, only the mean waiting time. Thus, the comparison is made over the mean waiting time for each station $2 \leq j \leq m$\footnote{Predicting the mean waiting time can be done by computing the average number of customers from \textbf{NN 2} and then deriving the mean waiting time via Littel's law. }.

To summarize, these are our main contributions:

\begin{itemize}
    \item Predicting the first departure moments and non-linear auto-correlation values from $GI/GI/1$ and $G/GI/1$ systems. 
    \item Predicting the steady-state number of customers in the system in a  $G/GI/1$.
    \item Conducting an empirical experiment that examines the effect of the inter-arrival and service moments and the inter-arrival auto-correlation on the steady-state distribution of a $G/GI/1$ queue.
\end{itemize}

The plan for the remainder of the paper is as follows. In Section~\ref{sec:lit_rev}, we provide a literature review. We formulate our problem in Section~\ref{sec:prob_formulation}.  We describe our training process in Section~\ref{sec:training_proc}. In Section~\ref{sec:experiments}, we detail our experiments and present their results in Section~\ref{sec:result}. We conclude the paper in Section~\ref{sec:conclusions}.

\section{Related Work}\label{sec:lit_rev}

Several papers have analyzed a tandem queueing system, some of which have a broader family of open queueing systems (OQNs) that can be reduced into a tandem queueing system.

\cite{Wang2024} study a $G/G/1$ queue in which inter-arrival times follow a non-renewal process. Their study is very similar to ours, except that it provides the waiting time moments while our model approximates the steady-state probabilities in each station at a tandem system. They also derive the moments and covariances
of the departure process. The authors use the MacLaurin series for their derivation. The main disadvantage of their method is its long running time, which mainly depends on the complexity of the (correlated) inter-arrival distribution.  The more complex the form of the inter-arrival time distribution and its auto-correlation, the longer the runtime. The calculation can be completed in a few seconds for simpler cases, such as an $M/M/1$ queue with the FGM copula. In more complex scenarios, such as with the Frank Copula, where inter-arrival times follow an Erlang distribution and service times are uniformly distributed, the runtime can extend to several days, typically longer than the simulation runtime.  

A similar study in \cite{Dai2022} considers a type of correlated queue in which the service time of a customer depends on the inter-arrival time between him and the previous customer. Using the Mckluren method, they compute both sojourn and inter-departure moments. These two studies are based on the results of~\cite{Gong_Hu_1992}, which was the first to use the Mckluren method for queueing analysis. Initially, it was used to derive the first moments of waiting time.

While those studies provide accurate predictions in a tandem system, they are extremely slow. As a result, approximation methods are required, which may not be as precise but can provide fast predictions. We distinguish between two types of studies: those that use ML-based approximations (like this study) and those that rely on queueing theory. We commence with the later type of papers. 

Tandem queueing systems are a special case of OQNs. Under the assumption of Poisson arrival processes and exponential service-time distributions, our tandem queueing network (in fact, any OQN) can be modeled as a Jackson network, which is relatively straightforward to analyze. This simplicity arises because the steady-state distribution of queue lengths follows a product form, meaning that the steady-state queue lengths are independent geometric random variables, similar to independent $M/M/1$ queues. 

Inspired by the product-form property of Markov OQNs, researchers have extensively studied decomposition approximations for non-Markov OQNs. This method breaks down the network into individual single-server queues, assuming the steady-state queue length processes are approximately independent. For instance, ~\cite{1094270} propose approximating each queue as a $GI/GI/1$ model, where the arrival and service processes are modeled as renewal processes, partially characterized by the mean and squared coefficient of variation (SCV, variance divided by the square of the mean) of interarrival or service times.

While decomposition approximations often yield good performance, it has been recognized that dependence on the arrival processes of internal flows can pose a significant challenge. The queueing analyzer (QNA) algorithm~\cite{6770716} includes an approximation to address this dependence. However, issues persisted, as demonstrated by comparisons between QNA and model simulations in studies such as~\cite{46511, doi:10.1287/mnsc.41.10.1704}. In~\cite{BARON2024106867}, the authors quantified under which settings ignoring the dependencies between arrivals incur substantial errors. \cite{CIVELEK20211031} checked the effect of linear correlation between inter-arrival and/or service time and waiting times. They use VARTA modeling to generate linear dependence between arrivals/services.  Then, they do simulations to examine the effect of the positive/negative correlations on the mean waiting time. Our study does not examine how it affects the mean waiting time or which direction but quantifies the effect of the $i^{th}$ lag correlation on the queue dynamics. As opposed to~\cite{CIVELEK20211031}, we do not consider only the linear correlations but higher order of polynomial correlations as well.  


Decomposition methods based on Markov arrival processes (MAPs) have been developed to tackle the dependence issue in arrival processes. MAPs, introduced by~\cite{Neuts_1979}, can capture dependencies among interarrival times (or service times) because a MAP is not a renewal process. \cite{HORVATH2010759} modeled each station using a $MAP/MAP/1$ framework, while~\cite{doi:10.1287/opre.1100.0893, Kim2011} approximated each queue with a $MMPP(2)/GI/1$ model, where the arrival process is a Markov-modulated Poisson process with two states, a specific type of MAP.

\cite{doi:10.1287/opre.2015.1367} introduced a novel robust queueing (RQ) approach for analyzing performance in single-server queues. The key idea behind RQ is to replace the traditional probabilistic framework with an appropriate uncertainty set, allowing for the analysis of worst-case performance in a deterministic manner. The authors used discrete-time Lindley's recursion to describe customer waiting times as the supremum of partial sums of interarrival and service times. The proposed uncertainty sets for these partial sums was derived using the central limit theorem and two-moment descriptions of the arrival and service processes.

Although the overall concept of RQ is straightforward and promising, challenges remain in identifying effective uncertainty sets and connecting the results back to the original queueing system. These challenges were addressed by~\cite{doi:10.1287/opre.2017.1649}, developing the robust queueing analyzer (RQNA). In that work, a new nonparametric RQ formulation was introduced to approximate the continuous-time workload process in a single-server queue, and it was proven that this approximation for the steady-state mean is asymptotically accurate in both light and heavy traffic conditions.

Other papers adopted an ML approach.  \cite{SAGRON2017161} suggested a regression-based variability function, which studied the departure process in a single server tandem queueing system via linear regression in a system with downtime events. Specifically, they approximated the departure variability of each station. Once the departure variability is approximated, they compute the downstream station's mean waiting time and the queue length. Their paper is designed for multi-class customers with different service rates. Their method predicts only the departure process variability while ignoring auto-correlation and higher moments. 

\cite{9869030} studied the mean waiting time of two station tandem queues, using only the first two moments of the inter-departure and service times. Instead of our model, they do not learn the departure process, instead,  study the waiting time of the second station directly. This implies their model cannot be extended to a more extensive tandem system.  

\cite{doi:10.1080/00207543.2021.1887536} used a GP regression to approximate an OQNs. Their model approximates the departure in the first two moments, first-lag auto-correlation, and the mean waiting time. Their model can handle a more general queueing system, as they allow superposition and splitting. The disadvantage of this paper is that they generated their inter-arrival times using only correlated Weibull sequences. By doing so, higher moments of the inter-arrival time are ignored; however, as the results of our model paper and of~\cite{sherzer23} demonstrate,  higher moments (i.e., beyond the first two) significantly impact the queue dynamics.

\section{Problem formulation}\label{sec:prob_formulation}
This section gives an exact formulation of the input and output for each NN. 
We have in total three neural networks, \textbf{NN 1}, \textbf{NN 2}, \textbf{NN 3}, which are utilized to predict the stationary behavior of an $m$ station tandem queueing system. 
For this task, we first present the following notation. 

\begin{itemize}
    \item $D_j$ is the stochastic stationary inter-departure from station $j \leq m$. 
    \item $S_j$ is the  stochastic stationary service time  at station $j \leq m$.
    \item $A$ is the stochastic stationary external inter-arrival to station 1.
    \item $m_{A}(i)$  is the true    $i^{th}$ moment of the  stationary  external inter-arrival.
    \item $m_{D,j}(i)$ ($\hat{m}_{D,j}(i)$) is the true (estimated)  $i^{th}$ moment of the  stationary inter-departure from  the $j^{th}$ station.
    \item $m_{S,j}(i)$ ($\hat{m}_{S,j}(i)$) is the true (estimated) $i^{th}$ moment of the stationary service time  of the $j^{th}$ station.
    \item $p_{l,j}$ ($\hat{p}_{l,j}$) is the true (estimated) stationary probability of having $l \geq 0$ customers at station $j$. 
\end{itemize}

The notation presented above is used to define the input and output for \textbf{NN 1}, \textbf{NN 2}, and \textbf{NN 3}, in Sections~\ref{sec:nn1},~\ref{sec:nn2}, and~\ref{sec:nn3}, respectively. We provide a summary in Section~\ref{sec:sum_prob_formulatiom}. 

\subsection{\textbf{NN 1}}\label{sec:nn1}

In \textbf{NN 1}, we approximate the departure process of a $GI/GI/1$ queue. The input is the first $n$ inter-arrival and service time moments. Formally, it is the series ($m_{A}(i),m_{S,1}(i)  $), for $i \leq n$. The output includes two parts. The first part is the first $n$ inter-departure estimated stationary moments: $\hat{m}_{D,1}(i)$, for $i \leq n$. The second part is auto-correlation values. For providing formal expression, we
 define the following: 

Let $D_{j,q}$ be the inter-departure time between the $q^{th}$ and the $(q + 1)^{th}$ customers under the $j^{th}$ station for some very large $q$\footnote{$q$ large enough so that the system is in steady-state. }. Let $\rho(a_1,a_2,k)$ be the 
$k^{th}$-lag auto-correlation between $D_{j,q}^{a_1}$ and $D_{j,q-k}^{a_2}$ for $a_1,a_2,k \in \mathbb{N}$. Formally, for $j \leq m$:

\begin{align} \label{eq:rho_definition}
    \rho(a_1,a_2,k, j) = \frac{Cov(D_{j,q}^{a_1}, D_{j,q-k}^{a_2})}{\sqrt{Var(D_{j,q}^{a_1})Var(D_{j,q-k}^{a_2})}}.
\end{align}

\noindent Let $\hat{\rho}(a_1,a_2,k, j)$ be its corresponding estimated value. 

As such, the second part of our output is the series $\hat{\rho}(a_1,a_2,k, 1)$ (as defined in Equation~\eqref{eq:rho_definition}), for  $ k \leq n_1 $,  and $a_1, a_2 \leq n_2$, for some $n_1, n_2 \in \mathbb{N}$.

\subsection{\textbf{NN 2}}\label{sec:nn2}

In this part, we wish to predict the stationary probabilities of the number of customers in the system. 
The input of such a system for station $2 \leq j \leq m$ is the inter-departure moments and auto-correlation from station $j-1$ and the first $n$ service time moments of station $j$. Formally, the input is $m_{D,j-1}(i)$ for $i \leq n$, $\rho(a_1,a_2,k, j-1)$, for  $ k \leq n_1 $,  and $a_1, a_2
 \leq n_2$, for some $n_1, n_2 \in \mathbb{N}$, and $m_{S,j}(i)$ for $i \leq n$. The output is $\hat{p}_{l,j}$ for $l\leq L$, where $L$ is some threshold in which we truncate the distribution such that the total probability of having more than $L$ customers is smaller than $\delta$. 
In our empirical evaluation, we used $\delta=10^{-3}$ to achieve the desired performance with $L=1500$ in all generated samples. This corresponds to covering the probability of between $0$ and $1499$ customers in the system; the last value of the output vector is set to the probability of having 1499 or \textit{more} customers.

\subsection{\textbf{NN 3}}\label{sec:nn3}

The input for \textbf{NN 3} with respect to station $j$ is the same as in \textbf{NN 2}, only the output changes. That is, the input is $m_{D,j-1}(i)$ for $i \leq n$, $\rho(a_1,a_2,k, j-1)$, for  $ k \leq n_1 $,  and $a_1, a_2 \leq n_2$, for some $n_1, n_2 \in \mathbb{N}$, and $m_{S,j}(i)$ for $i \leq n$. The output is: $\hat{m}_{D,j}(i)$ for $i \leq n$, $\hat{\rho}(a_1,a_2,k, j)$, for  $ k \leq n_1 $,  and $a_1, a_2 \leq n_2$, for some $n_1, n_2 \in \mathbb{N}$ and $2 \leq j \leq m$.

\subsection{Summary of problem formulation}\label{sec:sum_prob_formulatiom}

In Table~\ref{tab:sum_inp_out}, we summarize the input and output for each NN. 
\begin{table}[!htp]\centering
\caption{A summary of input-output for all NNs}\label{tab:sum_inp_out}
\scriptsize
\begin{tabular}{|l|r|r|r|r|r|}\toprule
NN 1 &Input &\multicolumn{3}{c|}{$m_{A}(i)$ and  $m_{S,1}(i)$ for $i \leq n$. } \\
NN 1 &Output &\multicolumn{3}{c|}{$\hat{m}_{D,1}(i)$ for $i \leq n$, $\hat{\rho}(a_1,a_2,k, 1)$, for  $ k \leq n_1 $,  and $a_1, a_2 \leq n_2$.} \\
\hline
NN 2 &Input &\multicolumn{3}{c|}{$m_{D,j-1}(i)$, $m_{S,j}(i)$ for $i \leq n$, $\rho(a_1,a_2,k, j-1)$, for  $ k \leq n_1 $, $a_1, a_2 \leq n_2$, and $2 \leq j \leq m$. } \\
NN 2 &Output &\multicolumn{3}{c|}{$\hat{p}_{l,j}$ for $l\leq L$.} \\
\hline
NN 3 &Input &\multicolumn{3}{c|}{$m_{D,j-1}(i)$, $m_{S,j}(i)$ for $i \leq n$, $\rho(a_1,a_2,k, j-1)$, for  $ k \leq n_1 $,  $a_1, a_2 \leq n_2$, and $2 \leq j \leq m$.} \\
NN 3 &Output &\multicolumn{3}{c|}{$\hat{m}_{D,j}(i)$ for $i \leq n$, $\hat{\rho}(a_1,a_2,k, j)$, for  $ k \leq n_1 $,  $a_1, a_2 \leq n_2$, and $2 \leq j \leq m$.} \\
\bottomrule
\end{tabular}
\end{table}

\section{Training process}\label{sec:training_proc}

This section consists of three parts. The first part details the training data generation process in Section~\ref{sec:data}. Then, we discuss the NN architecture in Section~\ref{sec:network}. For the third part, we present the loss function for all three of our NNs in Section~\ref{sec:loss_func}. This section is very important as it bridges between queueing theory and ML.

\subsection{Data generating}\label{sec:data}

In this subsection, we first describe the generation of input instances and then the labeling process. We then discuss the connection between the inter-arrival and service time distribution and the auto-correlation values. Finally, we will provide further information regarding the training, validation, and test sets used to calibrate our NNs.

\subsubsection{Input generation}\label{sec:input_generation}

As mentioned above, for training \textbf{NN 1}, \textbf{NN 2}, and \textbf{NN 3}, we are only required to simulate a two-station tandem system. In the tandem queueing system under consideration, the entire queueing dynamics are defined once the external inter-arrival distribution and the service times at each station are specified. The internal inter-arrival distributions follow directly from these assignments. Consequently, the task of data generation reduces to generating non-negative continuous time distributions. We wish to create a versatile distribution library that covers a wide range of distributions that are likely to include most real-world situations.

\begin{figure}
\centering
\includegraphics[scale = 0.35]{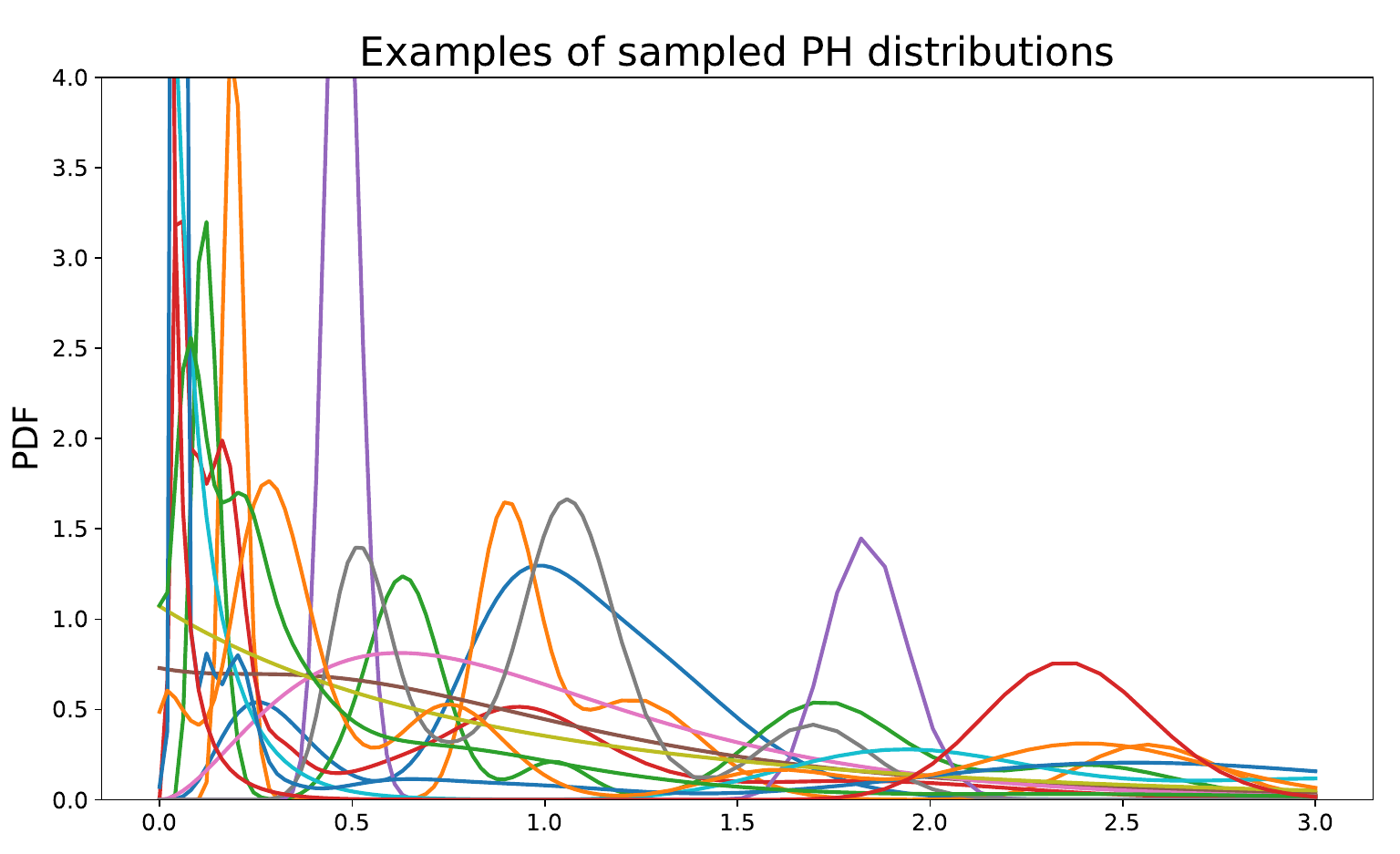}
\caption{An example of 15  distributions from~\cite{sherzer23}; the mean for all distributions is scaled to $1$.  }
\label{fig:pdf_dist}
\end{figure}

We employ the cutting-edge sampling method introduced in~\cite{sherzer23}. For distribution sampling, we use the Phase-Type (PH) family of distributions, which, as noted earlier, is dense within the class of all non-negative distributions. Although we sample from PH distributions up to a fixed degree (1000), this property enables us to generate a diverse training set in practice. In fact,~\cite{sherzer23} demonstrates broad distribution coverage using this approach, compares favorably with alternative sampling methods, and spans a wider family of real-life inter-arrival and service time distributions. Without the loss of generality, we scale the external inter-arrival distribution to be with a mean value of 1. The mean value of the service time is from 0 to 0.95, which means the utilization level for all stations is also under that range. The SCV values used for our training data range from 0.001 to 15. Although the generating method in~\cite{sherzer23} can generate larger SCV values, we limited the range as it is highly unlikely to get a larger SCV than 15 in real-life service systems, as shown in~\cite{sherzer23}. For illustration, we present 15 examples in Figure~\ref{fig:pdf_dist} of PH distribution to demonstrate the versatility of the PH generation procedure.

\subsubsection{Output generation}\label{sec:output_generation}

To label our data, we conduct simulations. For each instance, we run 60 million arrivals. To examine the accuracy of our simulations, we run 5000 tandem systems to explore the labeling consistency. We simulate each system ten times, each one including 60 million arrivals.  We then measure the length of a 95\%  confidence interval (CI) of selected output values for \textbf{NN 1}, \textbf{NN 2} and \textbf{NN 3}. For \textbf{NN 1} and \textbf{NN 3}, which predicts non-renewal departure processes, we examine the CI length of the first five moments and the first order linear auto-correlation (i.e., $a_1 = 1$, $a_2 = 1$, $k = 1$). For \textbf{NN 2}, we present CI lengths of the steady-state probabilities for 0 to 4 and the mean number of customers in the queue.

\begin{table}[!htp]\centering
\caption{Labeling CI interval}\label{tab:CI}
\scriptsize
\begin{tabular}{|l|r|r|r|r|r|r|r|}\toprule
&\multicolumn{5}{c}{Moments} &Utilization \\  \hline
Departures &1 &2 &3 &4 & 5 &(1,1,1) \\ \hline
NN 1 &0.0003 &0.0010 &0.0022 &0.0040 &0.0064 &0.0001 \\ \hline
NN 3 &0.0003 &0.0010 &0.0021 &0.0040 &0.0063 &0.0001 \\ \hline
&\multicolumn{5}{c}{Probabilites  } & Mean \\ \hline
Queue length &0 &1 &2 &3 &4 & \\ \hline
NN 2 &0.00017 &0.00005 &0.00004 &0.00003 &0.00003 &0.01973 \\ \hline
\bottomrule
\end{tabular}
\end{table}

The results are depicted in Table~\ref{tab:CI}. The CI interval lengths under all settings are low.  The mean CI interval length of the departure moments for \textbf{NN 1} and \textbf{NN 3} and similar and increase the moment order. The mean CI interval length of the mean number of customers is 0.01973, which reflects a relatively accurate labeling yet is slightly noisy. The average time a simulation with 60 million arrivals takes is 2.14 hours, using an Intel Xeon Gold 5115 Tray Processor with 128GB RAM. 

\begin{remark}\label{rem:log}
    We apply a log transformation for the moments to reduce the range of possible values. This aligns with the common practice in NN models to ensure that input values fall within the same range. Note that as we increase the order of moments, their values increase exponentially, thus even using standardized moments does not ensure range similarity (and adds to the computational burden). 
\end{remark}

\subsubsection{Auto-correlation space }\label{sec:auto-correlations}

The inter-departure auto-correlation values are not generated directly but result from the queueing model input, namely, external inter-arrival and service time distributions. Because auto-correlation values are critical in our analysis and are not generated straightforwardly, we show how they relate to the queueing model input. For this purpose, we examine the connection between the utilization level and the SCV values of the inter-arrival and service time distributions to the auto-correlation values. 

\begin{figure}
\centering
\includegraphics[scale=0.5]{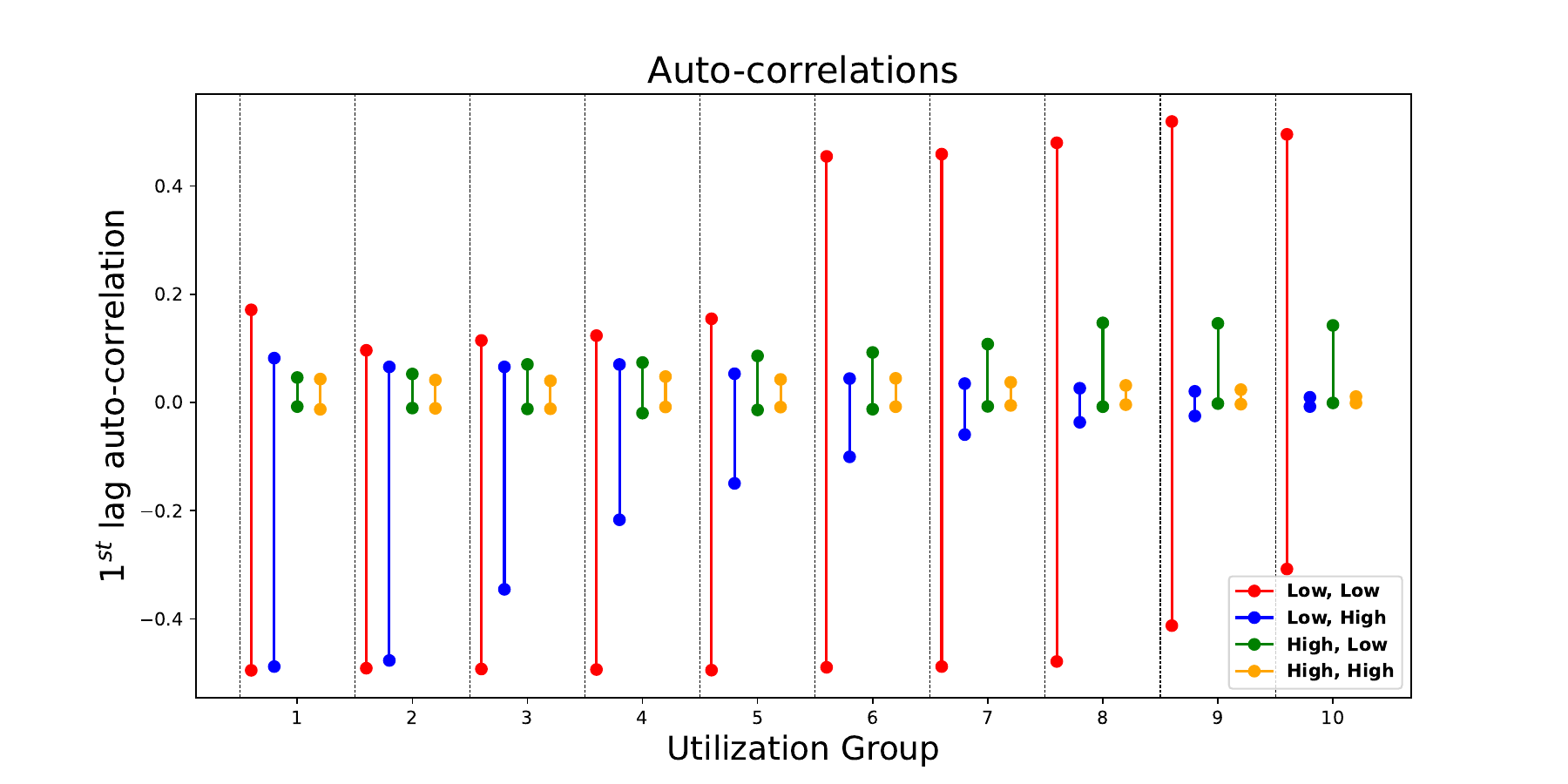}
\caption{Examining the $1^{st}$ lag linear auto-correlation range for different utilization and SCV groups. }
\label{fig:corrs_func}
\end{figure}

In Figure~\ref{fig:corrs_func}, we present the $1^{st}$ lag linear auto-correlation range for different utilization and SCV groups based on our training set. The utilization domain, ranging from 0 to 0.95, is divided into ten intervals. Each of the first nine intervals spans a width of 0.1, covering ranges [0, 0.1), [0.1, 0.2), ..., up to [0.8, 0.9). The final interval covers the range [0.9, 0.95]. The SCV values for inter-arrival and service time distributions are divided into two groups. The first contains lower SCV values ranging from 0 to 4. The second contains larger values, ranging from 4 to 15. We have $10\times 2 \times 2  = 40$ groups. 

The x-axis in Figure~\ref{fig:corrs_func} presents the different utilization groups. The utilization groups are ordered by their size, with the first group representing the range [0, 0.1) and the $10^{th}$ group representing [0.9,0.95].  Within each utilization group, we have four different SCV settings since we have two different distributions and two categories for each one. The first value is related to the SCV of the inter-arrival, and the second is the service time. For example, the red lines, which are under (Low, Low), represent the cases where both inter-arrival and service time distributions have SCV values that are lower than 4.

As depicted in the figure, auto-correlation values range from -0.5 to 0.5. We observe a strong effect of both utilization level and SCV values. To reach relatively large auto-correlation values that exceed 0.2, both inter-arrival and service time SCV need to be low, and the utilization level should be large, whereas, in our data, it should be larger than 0.5. 

Finally, we observe that the inter-arrival SCV value must be low to achieve low auto-correlation values. Meanwhile, for cases with high inter-arrival SCV values, auto-correlations are nearly 0.   
We do not aim to suggest why such an effect occurs; this is beyond the scope of this paper.

\subsubsection{Datasets: training, validation and test}\label{sec:datasets}

We use three datasets for each neural network: training, validation, and testing. The validation set optimizes the hyperparams such as learning rate, number of layers, etc.  We give the full specification of the hyperparams in Appendix~\ref{append:Fine-tuning}. 
The training, validation, and test set size is 300000, 10000, and 10000, respectively.

\subsection{Network Architecture}\label{sec:network}

We implement a fully connected feed-forward neural network (NN). In Appendix~\ref{append:architecture}, we provide a complete network specification. We use all NNs that are fully connected NN, often referred to as ANN (artificial neural network). The justification for selecting ANN instead of Convolutional Neural Networks (CNN) and Recurrent Neural Networks (RNN) is as follows.  CNNs are not ideal due to the small size of the input vector, and RNNs are more appropriate for time-dependent systems rather than stationary ones; see \cite{NIPS2015_b618c321, sherzer2024approximatinggtgi1queuesdeep} for further details.


\subsection{Loss function}\label{sec:loss_func}
In deep learning, training is typically performed in mini-batches, where a subset of the full training dataset (a "batch") is used to update the model parameters in each training iteration. We denote the batch size as $B$.  The batch size determines how many samples from the training data are processed before updating the model weights.   In \textbf{NN 1} and \textbf{NN 3}, which both model a departure, we use the same loss function and a separate one for \textbf{NN 2}. 
Recall that for training all NNs, we need only a two-station tandem system. As such, training \textbf{NN 1} is done via the departure of station 1,  training \textbf{NN 2} is done via the steady-state probabilities of station 2 and \textbf{NN 3} is done via the departure of station 2.

We introduce an additional index to include the $b^{th}$ instance in a batch. Let $m_{D,j,b}(i)$ ($\hat{m}_{D,j,b}(i)$) be the true (estimated)  $i^{th}$ moment of the  stationary inter-departure from  the $j^{th}$ station of the $b^{th}$ instance in a batch. For the notation of $\rho(a_1, a_2, k,j)$ and $\hat{\rho}(a_1, a_2, k,j)$, we add the index $b$ to include the instance within in a batch in the following way, $\rho_b(a_1, a_2, k,j)$ and $\hat{\rho_b}(a_1, a_2, k,j)$, respectively. The loss function of \textbf{NN 1} and \textbf{NN 3} is:

\begin{align*}
    &Loss_{ \{\textbf{NN 1, NN 3}\}}(\cdot) = \\ & \sum_{b=1}^B\! \left[ \! \frac{\alpha}{B} \! \left(\sum_{i=1}^n \left[log(m_{D_{j,b}}(i)) \!-\! log(\hat{m}_{D_{j,b}}(i))\right]^2 \right) \! + \!  \frac{(1\!-\!\alpha)}{B} \left(\sum_{k=1}^{n_1} \!\sum_{a_1=1}^{n_2}\!\sum_{a_2=1}^{n_2} \!\left[ \rho_b(a_1, a_2, k,j)\! - \!\hat{\rho_b}(a_1, a_2, k,j)\right]^2 \right) \right]
\end{align*}

\noindent where  $0 \leq \alpha \leq 1$, is a hyperparam and $j=1$ if we train \textbf{NN 1} and $j=2$ if we train \textbf{NN 3}.  Recall from Remark~\ref{rem:log} that we train over the moments log values, which explains why the log operator is in the loss function. This loss function is simply aggregating two parts: the moments and auto-correlation. For each part, we compute the MSE.

We next address the loss function of \textbf{NN 2}. Since in~\cite{sherzer23}, the authors also trained steady-state probabilities; we adopt the loss function used there in Equation (1). 
Let $p_{b,j,l}$ ($\hat{p}_{b,j,l}$) be the true (estimated) probability of having $l$ customers in the system, in the $b^{th}$ instance of the in the $j^{th}$ station.

We remind that for training \textbf{NN 2}, the stationary system distribution (output label), the distribution length is truncated at $L$. 
Our loss function for station $j = 2$ is given by,
\begin{equation}  \label{eq:loss}
     Loss_{\{\textbf{NN 2}\}}(\cdot) =   \frac{1}{B}\sum_{b=1}^{B}\sum_{l=0}^{L-1}|p_{b,j,l}-\hat{p}_{b,j,l}|+ 
     \frac{1}{B}\sum_{b=1}^{B}max_l(|p_{b,j,l}-
     \hat{p}_{b,j,l}|).
\end{equation}

\noindent See~\cite{sherzer23} for an explanation of the loss structure. 

\section{Experiments}\label{sec:experiments}

Our experiments have two goals: (i) examining the effect of the inter-depart and service time moments combined with auto-correlation values on queue dynamics and (ii) examining the accuracy of our method.  In Section~\ref{sec:acc_matircs}, we detail our accuracy matrix. Sections~\ref{sec:mom_auto_corr_anal} and~\ref{sec:perv_exp} discuss how goals (i) and (ii) are executed in our experiments, respectively. Section~\ref{sec:mod_comp_exp} details how our performances are compared to different approximation methods.

\subsection{Accuracy metrics}\label{sec:acc_matircs}

We distinguish between 3 types of measures that our model predicts. The  inter-departure moments (\textbf{NN 1}, \textbf{NN 3}), inter-departure auto-correlation values (\textbf{NN 1}, \textbf{NN 3}) and distribution PMF (\textbf{NN 2}). 

\noindent \textbf{Inter departure moments:} We use the  Mean Absolute Percentage Error (MAPE) for moments, which is the standard measure for such task~\cite{ https://doi.org/10.1002/nav.22010}. Formally,

\[
\text{MAPE} = \frac{1}{v} \sum_{i=1}^{v} \left| \frac{y_i - \hat{y}_i}{y_i} \right| \times 100,
\]

\noindent where $y_i$ and $\hat{y}_i$ and $v$ are the true and estimated value of the $i^{th}$ sample, respectively,  $v$ is the sample size.  

\noindent \textbf{Auto-correlation values:} they are typically small in magnitude. While MAPE is commonly used to evaluate prediction accuracy, it becomes less meaningful when applied to small values. Specifically, when the true values are close to zero, even small absolute differences between the predicted and true values can result in disproportionately large percentage errors. This can lead to misleading interpretations of model performance, overemphasizing relatively minor deviations.

Given this, we opt to use absolute error as a more appropriate measure for this task. The absolute error provides a direct measurement of the difference between the predicted and true auto-correlation values without normalizing by the magnitude of the true values. This avoids the inflation of error that occurs in percentage-based metrics when dealing with small values. By focusing on the absolute difference, we can better assess the model's ability to accurately predict auto-correlation without introducing the bias that percentage metrics introduce for small-scale data.

Therefore, we measure the model's performance using Mean Absolute Error (MAE), defined as:

\begin{align*}
    MAE = \frac{1}{v} \sum_{i=1}^{v} \left| y_i - \hat{y}_i\right|, 
\end{align*}

\noindent where $y_i$ and $\hat{y}_i$ and $v$ are the true and estimated value of the $i^{th}$ sample, respectively,  $v$ is the sample size.  

\noindent \textbf{Distribution PMF:} In this part, we seek a metric that effectively detects differences at any specified percentile, including the tails of the two distributions. We compute the MAPE of six percentile values: $25\%, 50\%, 75\%,90\%, 99\%, 99.9\%$. For presentation purposes, we present the $i^{th}$ estimated percentile of the number of customers in the system in the $j^{th}$ station by $\hat{F}^{-1}_{X_j}(i)$, where $X_j$ is the corresponding stochastic number of customers number of customers. 

We also evaluate the MAPE of the mean number of customers in the system. This allows us to assess the entire distribution accuracy via a single number. Further, as noted above, other approximation methods (e.g., ~\cite{doi:10.1287/opre.2015.1367} ,\cite{https://doi.org/10.1002/nav.22010}) compute only the mean value, which allows us to make comparisons. We denote the estimated mean number of customers in station $j$ by $m_{X_j}$. 

Finally, we have an additional metric we use in our experiments, which was initially introduced in~\cite{sherzer23}.  The Sum of Absolute Errors (SAE) between the predicted value and the labels. It is given by\begin{align}\label{eq:SAE}
  SAE =   \frac{1}{v}\sum_{i=1}^{v}\sum_{l=0}^{L-1}|p_{i,j,l}-\hat{p}_{i,j,l}|,
\end{align} where $v$ is the size of the dataset, and $j \in \{2,3,...,m\}$.   The $SAE$ is equivalent to the Wasserstein-1 measure and the first term of the loss function in Equation~\eqref{eq:loss}. The advantage of this metric is that it is highly sensitive to any difference between actual and predicted values of accuracy, and  it produces a single value, which simplifies the comparison of different models, making it particularly suitable for tuning the hyper-parameters of the model.  However, since this metric produces an absolute and not a relative value, it is more difficult to interpret when evaluating model accuracy. Thus, the only use of $SAE$ is for hyperparam tunning. As such, it is used for tunning $n$, $n_1$, and $n_2$, essential for the moments and auto-correlation analysis as presented in the following subsection.

\subsection{Moments and auto-correlation analysis}\label{sec:mom_auto_corr_anal}

In this section, we optimize the values of $n$, $n_1$, and $n_2$. Theoretically speaking, the larger those parameters are, the more information we have about the queueing system, which should increase accuracy. However, their information becomes less insignificant as these values increase. Furthermore, it is possible that increasing their values will be counterproductive as they insert more noise into the system. As such, we compare $SAE$ values while making predictions using \textbf{NN 2} under the test set.  We compare all possible combinations of $2  \leq n \leq 10$, and $0  \leq n_1, n_2 \leq 5$.

\subsection{Performence evaluation}\label{sec:perv_exp}

When evaluating our method for either inter-departure moments and auto-correlation values or steady-state probabilities, we wish to examine different utilization levels and different SCV of both the inter-arrival and service time distributions (see~\cite{https://doi.org/10.1002/nav.22010}).  For this purpose, we consider four utilization level intervals with equal lengths of 0.25\footnote{With one exception, the last interval is between 0.75 to 0.95}, and two groups of SCV, below and above 4. In total, we split into 16 different groups.

For examining the accuracy of \textbf{NN 1}, \textbf{NN 2}, and \textbf{NN 3}, we only need a 2-station tandem system. Whereas,  \textbf{NN 1} is examined via the station 1 departure process, \textbf{NN 2} is examined via station 2 stationary probabilities, and \textbf{NN 3} is examined via the station 2 departure process.  In Table~\ref{tab:experiment_discriptions}, we summarize all the experiments. Each experiment includes 16 combinations of the inter-arrive and service time SCV and utilization levels described above.

\begin{table}[!htp]\centering 
\caption{A summary of the accuracy exemperiments}\label{tab:experiment_discriptions}
\scriptsize
\begin{tabular}{lrrrrr}\toprule
Experiment \# &Network &Parameters &Metric &sample size \\
\hline
1 &NN 1 &$\hat{m}_{D,1}(i)$ for $i \leq n$ &MAPE &10000 \\
\hline
2 &NN 1 &$\hat{\rho}(a_1,a_2,k, 1)$, for $ k \leq n_1 $, and $a_1, a_2 \leq n_2$ &MAE &10000 \\
\hline
3 &NN 2 & $m_{X_2}$ and $F^{-1}_{p_j}(i)$ for $i \in \{25, 50,75,90,99,99.9\}$  &MAPE &10000 \\
\hline
4 &NN 3 &$\hat{m}_{D,1}(i)$ for $i \leq n$ &MAPE &10000 \\
\hline
5 &NN 3 &$\hat{\rho}(a_1,a_2,k, 2)$, for $ k \leq n_1 $, and $a_1, a_2 \leq n_2$ &MAE &10000 \\
\bottomrule
\end{tabular}
\end{table}

In experiment 3, as depicted in Table~\ref{tab:experiment_discriptions}, we examine the accuracy of our method for calculating the steady-state probabilities of a  $G/GI/1$ queue. As discussed above, we split the data into 16 groups of utilization levels and SCV values. We add four more scenarios, where we examine the accuracy of our method for different first-lag linear auto-correlation values (i.e., $\rho(1,1,1,1)$). As discussed in Section~\ref{sec:auto-correlations}, those values are highly affected by the utilization level and the SCV values. As such, we examine their accuracy values separately. The analysis in Section~\ref{sec:auto-correlations} shows that the first-lag linear auto-correlation values range is (-0.5,0.5). Accordingly, we split the data into four equal intervals of 0.25.

\subsection{Model comperison}\label{sec:mod_comp_exp}

In this section, we present a comparison between our model and the RQ model~\cite{doi:10.1287/opre.2015.1367} and RQNA~\cite{https://doi.org/10.1002/nav.22010}, which currently are state of the art. We replicate two comparisons in ~\cite{You19} on tandem queues. The first experiment is also a 2-station system, while the second is made on a 9-station system. For both settings, we use the following distributions. 

\begin{enumerate}
    \item Exponential ($M$) distribution with mean $1/\lambda$ and SCV = 1.
    \item Erlang ($E_k$) distribution with mean  $1/ \lambda$,  SCV$= 1/k$, i.e., the summation of k
i.i.d. exponential random variables, each with mean $1/(\lambda k)$.
    \item Hyperexponential $ (H_2(SCV))$  distribution, i.e., a mixture of two exponential distributions which can be parameterized by its first three moments or the mean $1/\lambda$, SCV  and
the ratio between the two components of the mean $r = p_1/\lambda_1 /(p_1/\lambda_1+p_2/\lambda_2)$ where
$\lambda_1 > \lambda_2$. Only the case $r = 0.5$ is considered.
\item  Log-normal ($LN(SCV)$ ) distribution with mean 1.
\item  Gamma ($G(4)$) distribution with mean 1 and SCV = $4$.
\end{enumerate}

\noindent \textbf{Comperison 1}: we look at various $GI_1/GI_2/1 \rightarrow \cdot /GI_3/1$ models, where $GI_1, GI_3 \in \{E_4, LN(4)\}$ and $GI_2 \in \{E_4, LN(0.25), M, H_2(4), LN(4), G(4)\}$. We denote the utilization levels of stations 1 and 2 by $\rho_1$ and $\rho_2$, respectively.  In this experiment, $\rho_1 \in  \{0.7, 0.9 \} $, and $\rho_2 \in \{0.11,0.16,...,, 0.91, 0.96\}$. In total, there are 864 scenarios. We aggregate the results into smaller groups.  Like the abovementioned experiments, we consider four utilization level $\rho_1$ intervals with equal lengths of 0.25. We keep the two values of $\rho_2$.  We consider two groups of SCV values: inter-arrive and service time. The first group will include distributions with SCV 0.25 and 1, which are low values, and the second group will include SCV values of 4, which are larger.
In total, we have 64 combinations.

\noindent \textbf{Comperison 2}
We compare our method's performance with the RQ model~\cite{doi:10.1287/opre.2015.1367} and RQNA~\cite{https://doi.org/10.1002/nav.22010}.  This example, introduced by~\cite{SURESH1990355}, serves as a benchmark.

Specifically, we examine an OQN with nine stations arranged in tandem, each having independent and identically distributed (i.i.d.) exponential service times. Station 1 has an external arrival process and a general renewal process with a rate of 1. The traffic intensities for the first eight queues are set at 0.6, while the last has a significantly higher traffic intensity of 0.9. As in~\cite{SURESH1990355}, two specific external renewal arrival processes are considered. The first is deterministic interarrival times with SCV 0, and the second is highly variable interarrival times with SCV 8, using a hyper-exponential distribution.

\section{Results}\label{sec:result}
In this section, we present the results of all our experiments. Section~\ref{sec:mom_anal} presents the moment analysis and auto-correlation results. In Section~\ref{sec:acc_res}, we present our model accuracy results for \textbf{NN 1}, \textbf{NN 2} and \textbf{NN 3}. We conclude this part in Section~\ref{sec:mod_comp}, where we compare our model against existing models in the literature. 

\subsection{Moment and auto-correlation analysis}\label{sec:mom_anal} 

We evaluate all possible combinations of $2 \leq n \leq 10$ and $0 \leq n_1, n_2 \leq 5$, where among all configurations of $n$, $n_1$, and $n_2$, the most accurate results were obtained with $n=5$, $n_1=2$, and $n_2=2$.
  
\begin{figure}[t!]
        \centering
        \begin{subfigure}[b]{0.45\textwidth}
            \centering
          \includegraphics[scale = 0.45]{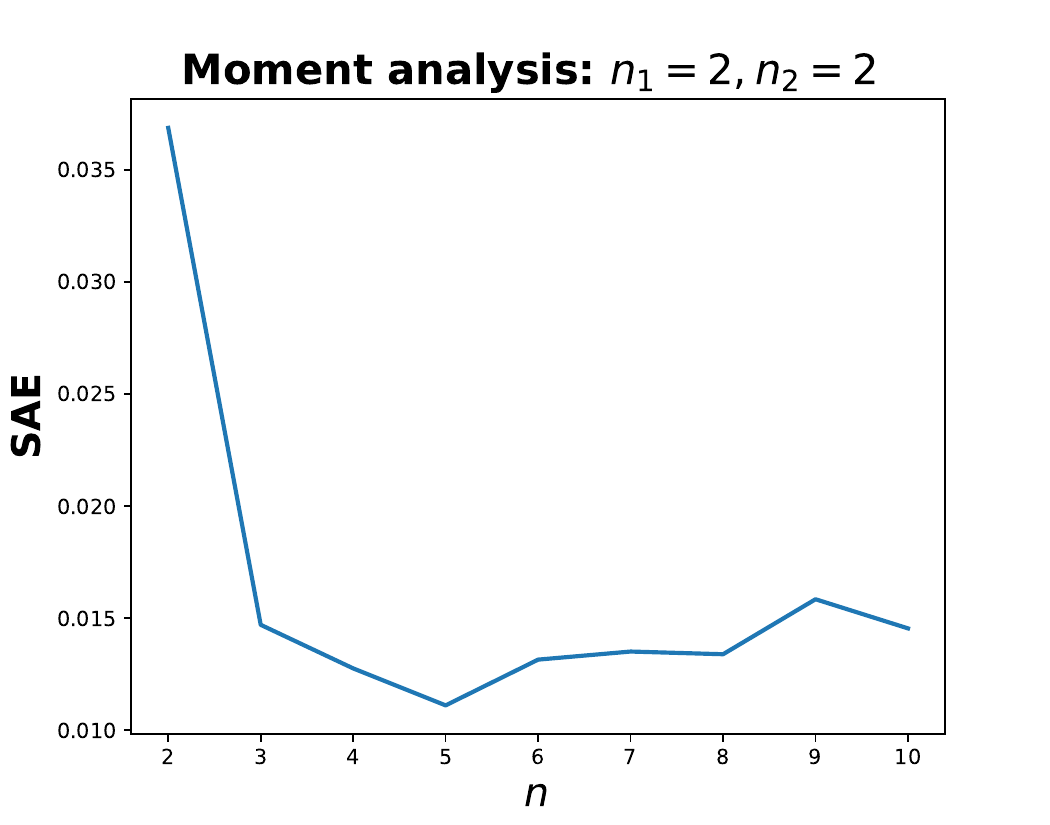}
            \caption[]%
            {{\small Moment analysis.}}    
            \label{fig:moms_anal}
        \end{subfigure}
        \hfill
        \begin{subfigure}[b]{0.45\textwidth}  
            \centering 
            \includegraphics[scale = 0.45]{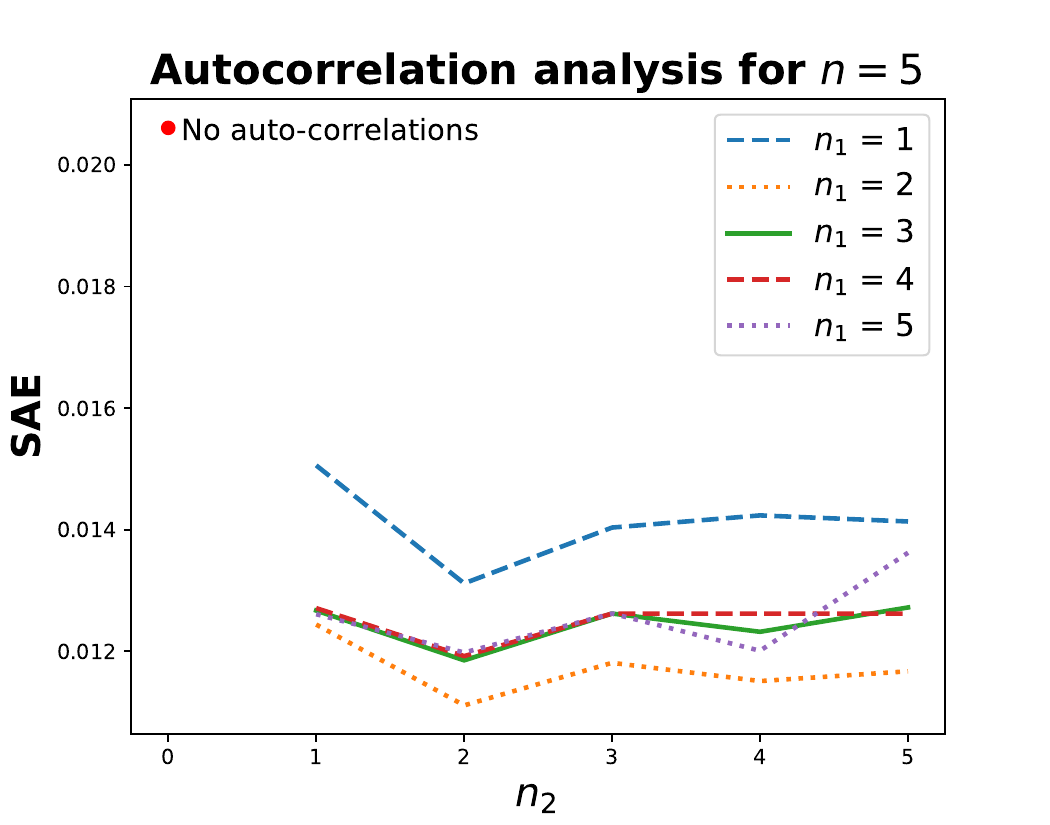}
            \caption[]%
            {{\small Auto-correlation analysis.}}    
            \label{fig:autocorr}
        \end{subfigure}
        \caption{Moment and auto-correlation analysis.}%
    \label{fig:mom_cor_anal}
    \end{figure}

The results depicted in Figure~\ref{fig:moms_anal} present the accuracy via the SAE measure as a function of the number of inter-arrival and service time moments $n$ that are used. The figure demonstrates the effect of each moment on the accuracy of a $G/GI/1$ queue. As the results show, more than five moments create more noise than added value. This supports previous results of a  $GI/GI/1$ queue, given in~\cite{sherzer23}. The SAE is not monotone from $n \geq 6$. We believe that this is because insignificant data, which is now inserted into the model as input, adds more noise than value. 

The results depicted in Figure~\ref{fig:autocorr} demonstrate the effect of the number of auto-correlation lags and their polynomial degree on the accuracy of a $G/GI/1$ queue. We observe that completely ignoring auto-correlation, as done in the red dot, creates a great loss of accuracy.
Interestingly, both $n_1$ and $n_2$ are optimized at 2. This means the 3-lag correlations and $3rd$ polynomial degree auto-correlation value have a negligible effect on the queue in that it inserts more noise than an added value.  The results of $n_1$ are supported by~\cite{CIVELEK20211031}. Their result, however, does not consider the $3^{rd}$ lag and beyond, nor non-linear correlations.

\begin{remark}
We do not conclude that the $6^{\text{th}}$ moment and higher, or the auto-correlation for $n_1, n_2 > 2$, do not affect the stationary queue dynamics. Rather, the lack of improvement in our model when including these suggests that, given our labeling accuracy and data volume, their marginal effect on the stationary queue dynamics cannot be expressed.
\end{remark}


\subsection{Accuracy results}\label{sec:acc_res}

In this section, we present our accuracy performance evaluation of $\textbf{NN 1}$, $\textbf{NN 2}$ and $\textbf{NN 3}$. We commence with the results of experiment 1, which are depicted in Table~\ref{tab:accuracy_results_depart_0_moms}. The LB and UB in Table~\ref{tab:accuracy_results_depart_0_moms}, represent the lower and upper bounds, respectively,  of the utilization levels and SCV values.

We observe an increase in the errors with the order of the moments. This may be partially due to the increase in the labeling error, as presented in Table~\ref{tab:CI}. There is no consistency of the error with respect to the utilization level. The inter-arrival SCV values seem more dominant than the service time SCV, where larger inter-arrival SCV comes with larger MAPE errors.

\begin{table}[!htp]\centering
\caption{Experiment 1: Accuracy results of \textbf{NN 1} moments}\label{tab:accuracy_results_depart_0_moms}
\scriptsize
\begin{tabular}{|c|c|c|c|c|c|c|c|c|c|c|}
\hline
\textbf{\#} & \multicolumn{2}{c|}{\textbf{Utilization}} & \multicolumn{2}{c|}{\textbf{SCV Arrival}} & \multicolumn{2}{c|}{\textbf{SCV services}} & \multicolumn{4}{c|}{\textbf{MAPE Departure Moments }} \\ 
\hline
 & \textbf{LB} & \textbf{UB} & \textbf{LB} & \textbf{UB} & \textbf{LB} & \textbf{UB} & \textbf{2} & \textbf{3} & \textbf{4} & \textbf{5} \\ 
\hline
1 & 0.00 & 0.25 & 0 & 4 & 0 & 4 & 0.79 & 1.08 & 1.54 & 1.81 \\ \hline
2 & 0.25 & 0.50 & 0 & 4 & 0 & 4 & 0.97 & 1.20 & 1.63 & 2.43 \\ \hline
3 & 0.50 & 0.75 & 0 & 4 & 0 & 4 & 0.70 & 1.24 & 1.87 & 2.74 \\ \hline
4 & 0.75 & 1.00 & 0 & 4 & 0 & 4 & 0.79 & 1.17 & 2.12 & 3.27 \\ \hline
5 & 0.00 & 0.25 & 0 & 4 & 4 & 15 & 1.45 & 1.59 & 2.14 & 3.83 \\ \hline
6 & 0.25 & 0.50 & 0 & 4 & 4 & 15 & 1.85 & 1.73 & 2.81 & 5.21 \\ \hline
7 & 0.50 & 0.75 & 0 & 4 & 4 & 15 & 1.42 & 2.12 & 4.58 & 7.44 \\ \hline
8 & 0.75 & 1.00 & 0 & 4 & 4 & 15 & 1.50 & 2.21 & 4.74 & 8.23 \\ \hline
9 & 0.00 & 0.25 & 4 & 15 & 0 & 4 & 3.26 & 2.62 & 3.80 & 5.21 \\ \hline
10 & 0.25 & 0.50 & 4 & 15 & 0 & 4 & 2.51 & 4.84 & 5.57 & 5.26 \\ \hline
11 & 0.50 & 0.75 & 4 & 15 & 0 & 4 & 3.37 & 5.30 & 5.01 & 7.29 \\ \hline
12 & 0.75 & 1.00 & 4 & 15 & 0 & 4 & 1.94 & 3.32 & 6.07 & 10.19 \\ \hline
13 & 0.00 & 0.25 & 4 & 15 & 4 & 15 & 3.84 & 2.72 & 3.98 & 6.25 \\ \hline
14 & 0.25 & 0.50 & 4 & 15 & 4 & 15 & 2.54 & 3.44 & 4.44 & 6.49 \\ \hline
15 & 0.50 & 0.75 & 4 & 15 & 4 & 15 & 1.46 & 2.73 & 5.95 & 8.42 \\ \hline
16 & 0.75 & 1.00 & 4 & 15 & 4 & 15 & 1.43 & 3.02 & 4.70 & 7.66 \\
\hline
\end{tabular}
\end{table}

Table~\ref{tab:auto_correlation_results_depart_0} depicts the MAE values for different auto-correlations. There is no substantial MAE difference between the first and second auto-correlation lags. Likewise, the second polynomial order is roughly as accurate as the first order. In General, the MAE values are very low.

\begin{table}[!htp]\centering
\caption{Experiment 2: accuracy results of \textbf{NN 1 }: Auto-correlations}\label{tab:auto_correlation_results_depart_0}
\scriptsize
\begin{tabular}{|p{0.5cm}|p{0.5cm}|p{0.5cm}|p{0.5cm}|p{0.5cm}|p{0.5cm}|p{0.5cm}|p{0.6cm}|p{0.6cm}|p{0.6cm}|p{0.6cm}|p{0.6cm}|p{0.6cm}|p{0.6cm}|p{0.6cm}|}
\hline
\textbf{\#} & \multicolumn{2}{|c|}{\textbf{Utilization}} & \multicolumn{2}{c|}{\textbf{SCV Arrival}} & \multicolumn{2}{c|}{\textbf{SCV services}} & \multicolumn{8}{c|}{\textbf{MAE Auto-correlation: $k, a_1, a_2$}} \\ 
\hline
 & \textbf{LB} & \textbf{UB} & \textbf{LB} & \textbf{UB} & \textbf{LB} & \textbf{UB} & \textbf{1,1,1} & \textbf{1,1,2} & \textbf{1,2,1} & \textbf{1,2,2} & \textbf{2,1,1} & \textbf{2,1,2} & \textbf{2,2,1} & \textbf{2,2,2} \\ 
\hline
1 & 0.00 & 0.25 & 0 & 4 & 0 & 4 & 0.006 & 0.006 & 0.006 & 0.005 & 0.003 & 0.002 & 0.003 & 0.002 \\ \hline
2 & 0.25 & 0.50 & 0 & 4 & 0 & 4 & 0.010 & 0.007 & 0.008 & 0.006 & 0.006 & 0.004 & 0.006 & 0.004 \\ \hline
3 & 0.50 & 0.75 & 0 & 4 & 0 & 4 & 0.010 & 0.011 & 0.010 & 0.008 & 0.006 & 0.006 & 0.006 & 0.006 \\ \hline
4 & 0.75 & 1.00 & 0 & 4 & 0 & 4 & 0.008 & 0.010 & 0.006 & 0.007 & 0.006 & 0.006 & 0.006 & 0.006 \\ \hline
5 & 0.00 & 0.25 & 0 & 4 & 4 & 15 & 0.016 & 0.009 & 0.011 & 0.007 & 0.011 & 0.004 & 0.009 & 0.005 \\ \hline
6 & 0.25 & 0.50 & 0 & 4 & 4 & 15 & 0.007 & 0.003 & 0.003 & 0.002 & 0.008 & 0.003 & 0.004 & 0.002 \\ \hline
7 & 0.50 & 0.75 & 0 & 4 & 4 & 15 & 0.003 & 0.001 & 0.001 & 0.001 & 0.002 & 0.001 & 0.001 & 0.001 \\ \hline
8 & 0.75 & 1.00 & 0 & 4 & 4 & 15 & 0.001 & 0.001 & 0.001 & 0.000 & 0.001 & 0.000 & 0.001 & 0.000 \\ \hline
9 & 0.00 & 0.25 & 4 & 15 & 0 & 4 & 0.002 & 0.001 & 0.001 & 0.001 & 0.001 & 0.001 & 0.001 & 0.001 \\ \hline
10 & 0.25 & 0.50 & 4 & 15 & 0 & 4 & 0.004 & 0.002 & 0.002 & 0.002 & 0.003 & 0.002 & 0.002 & 0.001 \\ \hline
11 & 0.50 & 0.75 & 4 & 15 & 0 & 4 & 0.002 & 0.001 & 0.001 & 0.001 & 0.002 & 0.001 & 0.001 & 0.001 \\ \hline
12 & 0.75 & 1.00 & 4 & 15 & 0 & 4 & 0.001 & 0.001 & 0.001 & 0.001 & 0.001 & 0.001 & 0.001 & 0.001 \\ \hline
13 & 0.00 & 0.25 & 4 & 15 & 4 & 15 & 0.002 & 0.002 & 0.001 & 0.002 & 0.002 & 0.001 & 0.001 & 0.001 \\ \hline
14 & 0.25 & 0.50 & 4 & 15 & 4 & 15 & 0.002 & 0.002 & 0.001 & 0.001 & 0.002 & 0.001 & 0.001 & 0.001 \\ \hline
15 & 0.50 & 0.75 & 4 & 15 & 4 & 15 & 0.001 & 0.001 & 0.001 & 0.001 & 0.001 & 0.001 & 0.001 & 0.001 \\ \hline
16 & 0.75 & 1.00 & 4 & 15 & 4 & 15 & 0.001 & 0.001 & 0.001 & 0.001 & 0.001 & 0.001 & 0.001 & 0.001 \\ \hline
\hline
\end{tabular}
\end{table}

The results of experiment 3 are presented in Tables~\ref{tab:steady_1} and Table~\ref{tab:nn2_cors}.   As the results in Table~\ref{tab:steady_1} suggest, the error is less than 5\% in most scenarios, even for a challenging task such as the 99.9\% percentile. We can see larger errors in the high percentiles for lower utility levels. We believe this is because when we have a small number of customers, even a slight deviation of a single customer may increase the MAPE value significantly.

Table~\ref{tab:nn2_cors} depicts the errors of our method as a function of the first lag linear auto-correlation values. The results of small auto-correlation values come with greater errors. This is likely because small auto-correlation values are associated with large SCV values of both the inter-arrival and service time distributions, as depicted in Figure~\ref{fig:corrs_func}.

\begin{table}[!htp]\centering
\caption{Experiment 3: Accuracy results for \textbf{NN 2}}\label{tab:steady_1}
\scriptsize
\begin{tabular}{|l|r|r|r|r|r|r|r|r|r|r|r|r|r|r|}\toprule
\textbf{} &\multicolumn{2}{|c|}{\textbf{Utilization}} &\multicolumn{2}{|c|}{\textbf{SCV Arrival}} &\multicolumn{2}{|c|}{\textbf{SCV service}} &\textbf{MAPE} &\multicolumn{6}{c|}{\textbf{MAPE Percentiles}} \\ \hline
\textbf{\#} &\textbf{LB} &\textbf{UB} &\textbf{LB} &\textbf{UB} &\textbf{LB} &\textbf{UB} &\textbf{Mean} &\textbf{25\%} &\textbf{50\%} &\textbf{75\%} &\textbf{90\%} &\textbf{99\%} &\textbf{99.9\%} \\ \hline \hline
1 &0 &0.25 &0 &4 &0 &4 &0.78 &0 &0 &0 &0.19 &1.67 &3.03 \\ \hline
2 &0.25 &0.5 &0 &4 &0 &4 &0.6 &0 &0 &0 &0 &1.2 &2.98 \\ \hline
3 &0 &0.25 &0 &4 &4 &15 &3.87 &0 &0 &0 &0.45 &3.08 &7.17 \\ \hline
4 &0.25 &0.5 &0 &4 &4 &15 &1.85 &0 &0 &0.52 &0.55 &2.36 &7.47 \\ \hline
5 &0 &0.25 &4 &15 &0 &4 &3.83 &0 &0 &0 &0.28 &4.33 &8.66 \\ \hline
6 &0.25 &0.5 &4 &15 &0 &4 &5.22 &0 &0 &1.75 &4.81 &7.46 &11.15 \\ \hline
7 &0 &0.25 &4 &15 &4 &15 &7.53 &0 &0 &0 &3.37 &5.85 &10.4 \\ \hline
8 &0.25 &0.5 &4 &15 &4 &15 &5.6 &0 &0 &3.4 &5.89 &6.5 &8.03 \\ \hline
9 &0.5 &0.75 &0 &4 &0 &4 &1.46 &0 &0.23 &1.02 &1.74 &1.29 &2.8 \\ \hline
10 &0.75 &1 &0 &4 &0 &4 &1.96 &0.12 &1.15 &1.34 &1.97 &2.74 &5.02 \\ \hline
11 &0.5 &0.75 &0 &4 &4 &15 &2.19 &0 &0 &2.3 &2.25 &2.76 &4.2 \\ \hline
12 &0.75 &1 &0 &4 &4 &15 &1.06 &2.13 &1.01 &1.42 &1.05 &1.83 &4.43 \\ \hline
13 &0.5 &0.75 &4 &15 &0 &4 &1.94 &0 &0.93 &1.76 &1.34 &2.54 &4.28 \\ \hline
14 &0.75 &1 &4 &15 &0 &4 &1.75 &0.96 &1.22 &1.48 &1.91 &2.53 &4.5 \\ \hline
15 &0.5 &0.75 &4 &15 &4 &15 &1.23 &0 &0 &0 &1.62 &1.27 &1.9 \\ \hline
16 &0.75 &1 &4 &15 &4 &15 &1.24 &0.59 &1.19 &1.45 &1.1 &1.72 &4.63 \\ \hline
\bottomrule
\end{tabular}
\end{table}

\begin{table}[!htp]\centering
\caption{Experiment 3: Accuracy results of \textbf{NN 2} by correlation}\label{tab:nn2_cors}
\scriptsize
\begin{tabular}{|l|r|r|r|r|r|r|r|r|r|r|}\toprule
&\multicolumn{2}{c}{\textbf{Correlation}} &\textbf{MAPE} &\multicolumn{6}{c|}{\textbf{MAPE Percentiles}} \\ \hline
\# &\textbf{LB} &\textbf{UB} &\textbf{Mean} &\textbf{25\%} &\textbf{50\%} &\textbf{75\%} 
&\textbf{90\%} &\textbf{99\%} &\textbf{99.9\%} \\ \hline
1 &-0.5 &-0.25 &0.02 &0.00 &0.87 &0.15 &1.29 &1.34 &3.34 \\  \hline
2 &-0.25 &0 &2.16 &0.77 &1.09 &1.36 &1.93 &2.96 &6.46 \\ \hline
3 &0 &0.25 &2.32 &1.27 &1.54 &1.78 &2.32 &3.81 &6.67 \\ \hline
4 &0.25 &0.5 &0.01 &0.00 &0.23 &0.04 &0.20 &1.35 &4.87 \\ \hline
\bottomrule
\end{tabular}
\end{table}

\begin{figure}
\centering
\includegraphics[scale=0.5]{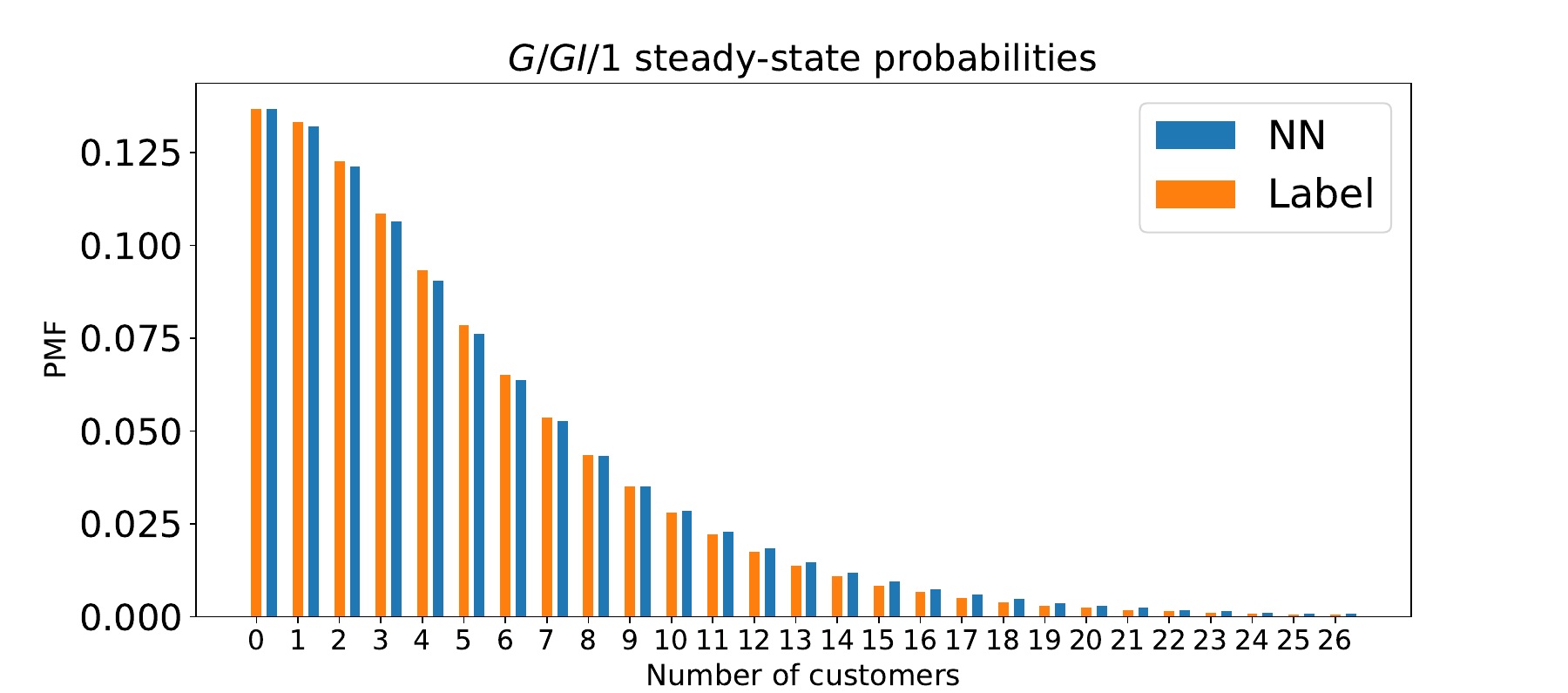}
\caption{Prediction example of a $G/GI/1$ queue. The orange bars are the labeled values given by simulations, and the blue bars are our NN predictions. }
\label{fig:gg1_exmaple}
\end{figure}

For illustration, we present an example of our steady-state probabilities prediction against the labeled values in Figure~\ref{fig:gg1_exmaple}. In this example, the MAPE of the mean value is 3.88, and the MAPE values for the 25\%, 50\%, 75\%, 90\%, 99\%, 99.9\% percentiles are 0,25,0,0 4.76 and 9.67, respectively. We deliberately choose a relatively non-accurate example to illustrate that even under these settings, our method appears to fit the true distribution well. 

The results for experiments 4 and 5 are presented in Table~\ref{tab:departure_moments_results} and \ref{tab:auto_correlation_results_depart_1}, where we examine the moment and auto-correlation accuracy, respectively. We observe the same results trends as we did in Tables~\ref{tab:accuracy_results_depart_0_moms} and ~\ref{tab:auto_correlation_results_depart_0}. As Table~\ref{tab:departure_moments_results} suggest, the MAPE errors are very small and under 3\% in most cases. Overall, the results under \textbf{NN 1} are slightly more accurate than under \textbf{NN 3}. This is sensible as the setting under \textbf{NN 3} is more complex due to the non-renewal arrival stream.

\begin{table}[!htp]\centering
\caption{Experiment 4: accuracy results of \textbf{NN 3 } moments}\label{tab:departure_moments_results}
\scriptsize
\begin{tabular}{|p{0.5cm}|p{0.5cm}|p{0.5cm}|p{0.5cm}|p{0.5cm}|p{0.5cm}|p{0.5cm}|p{0.6cm}|p{0.6cm}|p{0.6cm}|p{0.6cm}|}
\hline
\textbf{\#} & \multicolumn{2}{c|}{\textbf{Utilization}} & \multicolumn{2}{c|}{\textbf{SCV Arrival}} & \multicolumn{2}{c|}{\textbf{SCV services}} & \multicolumn{4}{c|}{\textbf{Departure Moments}} \\ 
\hline
 & \textbf{LB} & \textbf{UB} & \textbf{LB} & \textbf{UB} & \textbf{LB} & \textbf{UB} & \textbf{2} & \textbf{3} & \textbf{4} & \textbf{5} \\ 
\hline
1 & 0.00 & 0.25 & 0 & 4 & 0 & 4 & 1.53 & 1.46 & 1.66 & 2.12 \\ \hline
2 & 0.25 & 0.50 & 0 & 4 & 0 & 4 & 1.99 & 2.27 & 2.56 & 3.35 \\ \hline
3 & 0.50 & 0.75 & 0 & 4 & 0 & 4 & 2.06 & 2.20 & 2.74 & 3.45 \\ \hline
4 & 0.75 & 1.00 & 0 & 4 & 0 & 4 & 1.65 & 2.09 & 2.97 & 4.57 \\ \hline
5 & 0.00 & 0.25 & 0 & 4 & 4 & 15 & 2.63 & 2.57 & 3.48 & 4.02 \\ \hline
6 & 0.25 & 0.50 & 0 & 4 & 4 & 15 & 2.94 & 3.29 & 3.21 & 5.19 \\ \hline
7 & 0.50 & 0.75 & 0 & 4 & 4 & 15 & 2.35 & 2.70 & 3.77 & 6.29 \\ \hline
8 & 0.75 & 1.00 & 0 & 4 & 4 & 15 & 1.77 & 2.27 & 3.36 & 5.58 \\ \hline
9 & 0.00 & 0.25 & 4 & 15 & 0 & 4 & 2.32 & 1.82 & 1.97 & 2.13 \\ \hline
10 & 0.25 & 0.50 & 4 & 15 & 0 & 4 & 2.95 & 2.59 & 2.75 & 3.14 \\ \hline
11 & 0.50 & 0.75 & 4 & 15 & 0 & 4 & 3.22 & 3.34 & 3.72 & 4.46 \\ \hline
12 & 0.75 & 1.00 & 4 & 15 & 0 & 4 & 2.26 & 2.84 & 3.80 & 5.47 \\ \hline
13 & 0.00 & 0.25 & 4 & 15 & 4 & 15 & 2.21 & 1.42 & 1.71 & 2.72 \\ \hline
14 & 0.25 & 0.50 & 4 & 15 & 4 & 15 & 2.70 & 3.28 & 2.77 & 3.27 \\ \hline
15 & 0.50 & 0.75 & 4 & 15 & 4 & 15 & 1.92 & 2.62 & 4.52 & 6.99 \\ \hline
16 & 0.75 & 1.00 & 4 & 15 & 4 & 15 & 1.77 & 2.10 & 3.02 & 5.02 \\ \hline
\hline
\end{tabular}
\end{table}

\begin{table}[!htp]\centering
\caption{Experiment 5: accuracy results of \textbf{NN 3 } Auto-correlations}\label{tab:auto_correlation_results_depart_1}
\scriptsize
\begin{tabular}{|p{0.5cm}|p{0.5cm}|p{0.5cm}|p{0.5cm}|p{0.5cm}|p{0.5cm}|p{0.5cm}|p{0.6cm}|p{0.6cm}|p{0.6cm}|p{0.6cm}|p{0.6cm}|p{0.6cm}|p{0.6cm}|p{0.6cm}|}
\hline
\textbf{\#} & \multicolumn{2}{c|}{\textbf{Utilization}} & \multicolumn{2}{c|}{\textbf{SCV Arrival}} & \multicolumn{2}{c|}{\textbf{SCV services}} & \multicolumn{8}{c|}{\textbf{Auto-correlation: $k, a_1, a_2$}} \\ 
\hline
 & \textbf{LB} & \textbf{UB} & \textbf{LB} & \textbf{UB} & \textbf{LB} & \textbf{UB} & \textbf{1,1,1} & \textbf{1,1,2} & \textbf{1,2,1} & \textbf{1,2,2} & \textbf{2,1,1} & \textbf{2,1,2} & \textbf{2,2,1} & \textbf{2,2,2} \\ 
\hline
1 & 0.00 & 0.25 & 0 & 4 & 0 & 4 & 0.010 & 0.013 & 0.013 & 0.010 & 0.006 & 0.005 & 0.007 & 0.004 \\ \hline
2 & 0.25 & 0.50 & 0 & 4 & 0 & 4 & 0.013 & 0.013 & 0.010 & 0.008 & 0.011 & 0.007 & 0.011 & 0.007 \\ \hline
3 & 0.50 & 0.75 & 0 & 4 & 0 & 4 & 0.008 & 0.005 & 0.004 & 0.003 & 0.005 & 0.003 & 0.004 & 0.003 \\ \hline
4 & 0.75 & 1.00 & 0 & 4 & 0 & 4 & 0.005 & 0.004 & 0.004 & 0.003 & 0.003 & 0.002 & 0.003 & 0.002 \\ \hline
5 & 0.00 & 0.25 & 0 & 4 & 4 & 15 & 0.016 & 0.014 & 0.008 & 0.010 & 0.012 & 0.007 & 0.012 & 0.006 \\ \hline
6 & 0.25 & 0.50 & 0 & 4 & 4 & 15 & 0.010 & 0.005 & 0.005 & 0.003 & 0.007 & 0.004 & 0.006 & 0.002 \\ \hline
7 & 0.50 & 0.75 & 0 & 4 & 4 & 15 & 0.003 & 0.001 & 0.001 & 0.001 & 0.002 & 0.001 & 0.002 & 0.001 \\ \hline
8 & 0.75 & 1.00 & 0 & 4 & 4 & 15 & 0.001 & 0.001 & 0.001 & 0.001 & 0.001 & 0.001 & 0.001 & 0.000 \\ \hline
9 & 0.00 & 0.25 & 4 & 15 & 0 & 4 & 0.004 & 0.002 & 0.002 & 0.002 & 0.003 & 0.002 & 0.002 & 0.002 \\ \hline
10 & 0.25 & 0.50 & 4 & 15 & 0 & 4 & 0.004 & 0.002 & 0.002 & 0.002 & 0.003 & 0.002 & 0.001 & 0.001 \\ \hline
11 & 0.50 & 0.75 & 4 & 15 & 0 & 4 & 0.007 & 0.001 & 0.004 & 0.001 & 0.004 & 0.001 & 0.001 & 0.002 \\ \hline
12 & 0.75 & 1.00 & 4 & 15 & 0 & 4 & 0.002 & 0.001 & 0.001 & 0.001 & 0.001 & 0.001 & 0.001 & 0.001 \\ \hline
13 & 0.00 & 0.25 & 4 & 15 & 4 & 15 & 0.004 & 0.003 & 0.002 & 0.002 & 0.003 & 0.003 & 0.001 & 0.002 \\ \hline
14 & 0.25 & 0.50 & 4 & 15 & 4 & 15 & 0.005 & 0.002 & 0.002 & 0.001 & 0.004 & 0.002 & 0.002 & 0.001 \\ \hline
15 & 0.50 & 0.75 & 4 & 15 & 4 & 15 & 0.002 & 0.002 & 0.002 & 0.001 & 0.002 & 0.002 & 0.002 & 0.001 \\ \hline
16 & 0.75 & 1.00 & 4 & 15 & 4 & 15 & 0.001 & 0.001 & 0.001 & 0.001 & 0.001 & 0.001 & 0.001 & 0.001 \\ \hline
\hline
\end{tabular}
\end{table}

\subsection{Model comperison}\label{sec:mod_comp}

Tables~\ref{tab:2_comp_low} and~\ref{tab:2_comp_high} depicts \textbf{Comperison 1}, where Table~\ref{tab:2_comp_low} refers to cases where $\rho_1  \leq  0.5$ and Table~\ref{tab:2_comp_high} refers to cases where $\rho_1  \geq  0.5$. The values under the columns NN, RQNA, and RQ represent the MAPE value of the mean sojourn time of all three models. In each row, the model with the smallest MAPE is highlighted in red. Out of 64 scenarios, our NN model has the lowest MAPE 50 times. Furthermore, our model's worst MAPE value is 8.23, while in the RQNA and RQ it is 26.55 and 45.57, respectively. 

\begin{table}[!htp]\centering
\caption{2 Station comparison: low utility levels}\label{tab: }\label{tab:2_comp_low}
\scriptsize
\begin{tabular}{|l|r|r|r|r|r|r|r|r|r|r|}\toprule
&\textbf{Utilization} &\multicolumn{2}{c|}{\textbf{Utilization Station 2}} &\multicolumn{3}{c|}{\textbf{SCV}} &\multicolumn{3}{c|}{\textbf{MAPE}} \\ \hline
\textbf{\#} &\textbf{Station 1} &\textbf{ LB} &\textbf{UB} &\textbf{Arrival} &\textbf{Service 1} &\textbf{Service 2} &\textbf{NN} &\textbf{RQNA} &\textbf{RQ} \\ \hline
\textbf{1} &0.7 &0 &0.25 &Low &Low &Low &\textcolor{red}{\textbf{2.07}} &12.01 &11.92 \\ \hline
\textbf{2} &0.7 &0 &0.25 &High &Low &Low &\textcolor{red}{\textbf{3.80}} &11.40 &11.42 \\ \hline
\textbf{3} &0.7 &0 &0.25 &Low &Low &High &\textcolor{red}{\textbf{2.14}} &4.00 &3.90 \\ \hline
\textbf{4} &0.7 &0 &0.25 &High &Low &High &2.53 &2.60 &\textcolor{red}{\textbf{2.42}} \\ \hline
\textbf{5} &0.7 &0 &0.25 &Low &High &Low &\textcolor{red}{\textbf{2.04}} &6.59 &10.46 \\ \hline
\textbf{6} &0.7 &0 &0.25 &High &High &Low &\textcolor{red}{\textbf{1.15}} &16.47 &17.69 \\ \hline
\textbf{7} &0.7 &0 &0.25 &Low &High &High &2.88 &\textcolor{red}{\textbf{2.77}} &3.70 \\ \hline
\textbf{8} &0.7 &0 &0.25 &High &High &High &\textcolor{red}{\textbf{1.55}} &8.30 &8.48 \\ \hline
\textbf{9} &0.9 &0 &0.25 &Low &Low &Low &\textcolor{red}{\textbf{2.87}} &11.34 &11.34 \\ \hline
\textbf{10} &0.9 &0 &0.25 &High &Low &Low &\textcolor{red}{\textbf{4.13}} &11.06 &10.68 \\ \hline
\textbf{11} &0.9 &0 &0.25 &Low &Low &High &\textcolor{red}{\textbf{2.82}} &4.15 &4.17 \\ \hline
\textbf{12} &0.9 &0 &0.25 &High &Low &High &4.19 &3.68 &\textcolor{red}{\textbf{3.63}} \\ \hline
\textbf{13} &0.9 &0 &0.25 &Low &High &Low &\textcolor{red}{\textbf{5.08}} &12.78 &14.00 \\ \hline
\textbf{14} &0.9 &0 &0.25 &High &High &Low &\textcolor{red}{\textbf{7.20}} &15.74 &16.28 \\ \hline
\textbf{15} &0.9 &0 &0.25 &Low &High &High &\textcolor{red}{\textbf{4.73}} &5.63 &6.20 \\ \hline
\textbf{16} &0.9 &0 &0.25 &High &High &High &\textcolor{red}{\textbf{6.65}} &7.88 &7.91 \\ \hline
\textbf{17} &0.7 &0.25 &0.5 &Low &Low &Low &\textcolor{red}{\textbf{2.81}}&18.45 &18.21 \\ \hline
\textbf{18} &0.7 &0.25 &0.5 &High &Low &Low &\textcolor{red}{\textbf{6.39}} &14.32 &12.69 \\ \hline
\textbf{19} &0.7 &0.25 &0.5 &Low &Low &High &2.91 &1.92 &\textcolor{red}{\textbf{1.87}} \\  \hline
\textbf{20} &0.7 &0.25 &0.5 &High &Low &High &\textcolor{red}{\textbf{4.24} }&5.87 &4.85 \\  \hline
\textbf{21} &0.7 &0.25 &0.5 &Low &High &Low &\textcolor{red}{\textbf{3.18} }&12.67 &13.58 \\ \hline
\textbf{22} &0.7 &0.25 &0.5 &High &High &Low &\textcolor{red}{\textbf{2.12}} &22.69 &20.86 \\ \hline
\textbf{23} &0.7 &0.25 &0.5 &Low &High &High &3.35 &4.44 &\textcolor{red}{\textbf{3.02}} \\ \hline
\textbf{24} &0.7 &0.25 &0.5 &High &High &High &\textcolor{red}{\textbf{1.86} }&8.44 &6.88 \\ \hline
\textbf{25} &0.9 &0.25 &0.5 &Low &Low &Low &\textcolor{red}{\textbf{4.79} }&18.37 &18.23 \\ \hline
\textbf{26} &0.9 &0.25 &0.5 &High &Low &Low &\textcolor{red}{\textbf{7.50} }&17.78 &17.02 \\ \hline
\textbf{27} &0.9 &0.25 &0.5 &Low &Low &High &3.05 &1.46 &\textcolor{red}{\textbf{1.44}} \\ \hline
\textbf{28} &0.9 &0.25 &0.5 &High &Low &High &5.07 &2.23 &\textcolor{red}{\textbf{1.84}} \\ \hline
\textbf{29} &0.9 &0.25 &0.5 &Low &High &Low &\textcolor{red}{\textbf{5.81}} &14.16 &16.44 \\ \hline
\textbf{30} &0.9 &0.25 &0.5 &High &High &Low &\textcolor{red}{\textbf{8.23}} &19.46 &18.66 \\ \hline
\textbf{31} &0.9 &0.25 &0.5 &Low &High &High &\textcolor{red}{\textbf{3.84}} &4.60 &5.08 \\ \hline
\textbf{32} &0.9 &0.25 &0.5 &High &High &High &\textcolor{red}{\textbf{4.94}}&7.38 &6.71 \\ \hline
\bottomrule
\end{tabular}
\end{table}

\begin{table}[!htp]\centering

\caption{2 Station comparison: high utility levels}\label{tab:2_comp_high}
\scriptsize
\begin{tabular}{|l|r|r|r|r|r|r|r|r|r|r|}\toprule 
&\textbf{Utilization} &\multicolumn{2}{c|}{\textbf{Utilization Station 2}} &\multicolumn{3}{c|}{\textbf{SCV}} &\multicolumn{3}{c|}{\textbf{MAPE}} \\ \hline
\textbf{\#} &\textbf{Station 1} &\textbf{ LB} &\textbf{UB} &\textbf{Arrival} &\textbf{Service 1} &\textbf{Service 2} &\textbf{NN} &\textbf{RQNA} &\textbf{RQ} \\ \hline
\textbf{1} &0.7 &0.5 &0.75 &Low &Low &Low &\textcolor{red}{\textbf{2.06}} &7.02 &7.41 \\ \hline
\textbf{2} &0.7 &0.5 &0.75 &High &Low &Low &\textcolor{red}{\textbf{3.84}} &9.27 &11.76 \\ \hline
\textbf{3} &0.7 &0.5 &0.75 &Low &Low &High &1.44 &1.98 &\textcolor{red}{\textbf{1.44}} \\ \hline
\textbf{4} &0.7 &0.5 &0.75 &High &Low &High &\textcolor{red}{\textbf{1.69}} &4.24 &3.55 \\ \hline
\textbf{5} &0.7 &0.5 &0.75 &Low &High &Low &\textcolor{red}{\textbf{1.98}} &10.49 &6.24 \\ \hline
\textbf{6} &0.7 &0.5 &0.75 &High &High &Low &\textcolor{red}{\textbf{2.08}} &9.02 &6.29 \\ \hline
\textbf{7} &0.7 &0.5 &0.75 &Low &High &High &\textcolor{red}{\textbf{1.19}} &2.82 &2.16 \\ \hline
\textbf{8} &0.7 &0.5 &0.75 &High &High &High &\textcolor{red}{\textbf{0.68}} &2.27 &3.11 \\ \hline
\textbf{9} &0.9 &0.5 &0.75 &Low &Low &Low &\textcolor{red}{\textbf{3.65}} &7.11 &7.26 \\ \hline
\textbf{10} &0.9 &0.5 &0.75 &High &Low &Low &5.24 &7.51 &\textcolor{red}{\textbf{4.94}} \\ \hline
\textbf{11} &0.9 &0.5 &0.75 &Low &Low &High &1.51 &1.36 &\textcolor{red}{\textbf{1.13}} \\ \hline
\textbf{12} &0.9 &0.5 &0.75 &High &Low &High &\textcolor{red}{\textbf{1.65}} &3.45 &1.67 \\ \hline
\textbf{13} &0.9 &0.5 &0.75 &Low &High &Low &\textcolor{red}{\textbf{3.04}} &6.64 &5.92 \\ \hline
\textbf{14} &0.9 &0.5 &0.75 &High &High &Low &\textcolor{red}{\textbf{4.62}} &5.64 &7.16 \\ \hline
\textbf{15} &0.9 &0.5 &0.75 &Low &High &High &\textcolor{red}{\textbf{1.50}} &3.57 &3.25 \\ \hline
\textbf{16} &0.9 &0.5 &0.75 &High &High &High &\textcolor{red}{\textbf{1.76}} &1.84 &3.34 \\ \hline
\textbf{17} &0.7 &0.75 &1 &Low &Low &Low &\textcolor{red}{\textbf{2.75}} &8.94 &6.17 \\ \hline
\textbf{18} &0.7 &0.75 &1 &High &Low &Low &\textcolor{red}{\textbf{2.46}} &25.04 &33.65 \\ \hline
\textbf{19} &0.7 &0.75 &1 &Low &Low &High &3.19 &1.27 &\textcolor{red}{\textbf{1.16}} \\ \hline
\textbf{20} &0.7 &0.75 &1 &High &Low &High &\textcolor{red}{\textbf{2.69}} &7.71 &9.62 \\ \hline
\textbf{21} &0.7 &0.75 &1 &Low &High &Low &\textcolor{red}{\textbf{3.42}} &15.69 &10.04 \\ \hline
\textbf{22} &0.7 &0.75 &1 &High &High &Low &\textcolor{red}{\textbf{3.23}} &4.69 &11.08 \\ \hline
\textbf{23} &0.7 &0.75 &1 &Low &High &High &\textcolor{red}{\textbf{3.25}} &6.35 &4.52 \\ \hline
\textbf{24} &0.7 &0.75 &1 &High &High &High &3.64 &\textcolor{red}{\textbf{2.31}} &4.88 \\ \hline
\textbf{25} &0.9 &0.75 &1 &Low &Low &Low &\textcolor{red}{\textbf{3.75}} &6.27 &4.59 \\ \hline
\textbf{26} &0.9 &0.75 &1 &High &Low &Low &\textcolor{red}{\textbf{6.24}} &26.55 &45.57 \\ \hline
\textbf{27} &0.9 &0.75 &1 &Low &Low &High &3.74 &1.02 &\textcolor{red}{\textbf{0.90}} \\ \hline
\textbf{28} &0.9 &0.75 &1 &High &Low &High &\textcolor{red}{\textbf{4.60}} &4.91 &8.62 \\ \hline
\textbf{29} &0.9 &0.75 &1 &Low &High &Low &\textcolor{red}{\textbf{2.49}} &7.66 &7.26 \\ \hline
\textbf{30} &0.9 &0.75 &1 &High &High &Low &\textcolor{red}{\textbf{3.80}} &4.52 &11.68 \\ \hline
\textbf{31} &0.9 &0.75 &1 &Low &High &High &\textcolor{red}{\textbf{2.98}} &4.36 &4.17 \\ \hline
\textbf{32} &0.9 &0.75 &1 &High &High &High &3.91 &\textcolor{red}{\textbf{1.55}} &4.98 \\ \hline
\bottomrule
\end{tabular}
\end{table}

\begin{figure}[t!]
        \centering
        \begin{subfigure}[b]{0.45\textwidth}
            \centering
          \includegraphics[scale = 0.45]{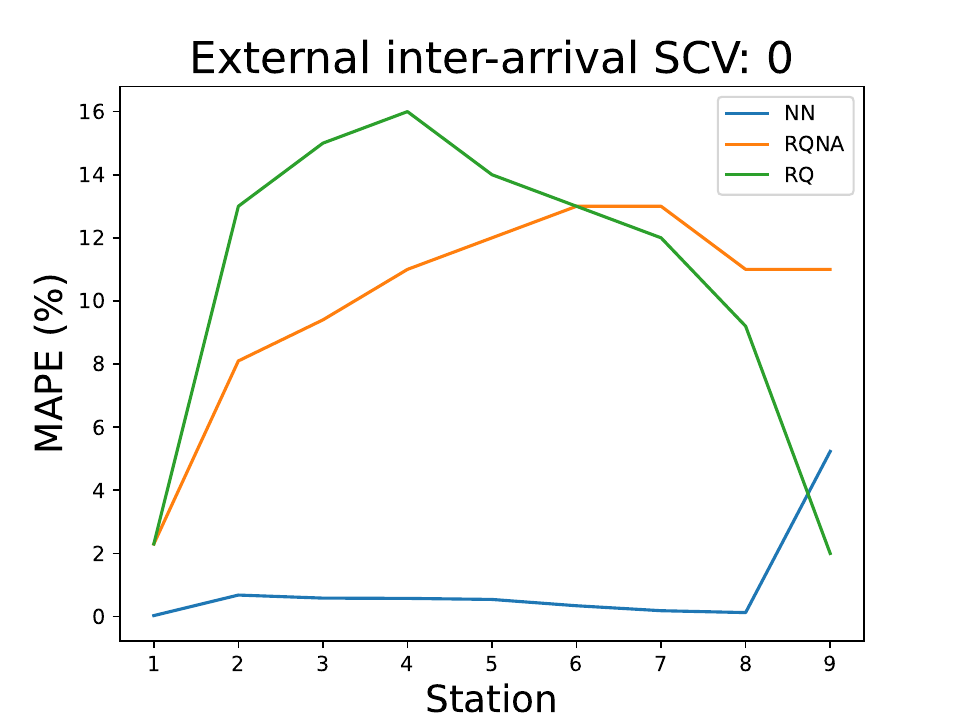}
            \caption[]%
            {{\small Determenstic arrival.}}    
            \label{fig:deter}
        \end{subfigure}
        \hfill
        \begin{subfigure}[b]{0.45\textwidth}  
            \centering 
            \includegraphics[scale = 0.45]{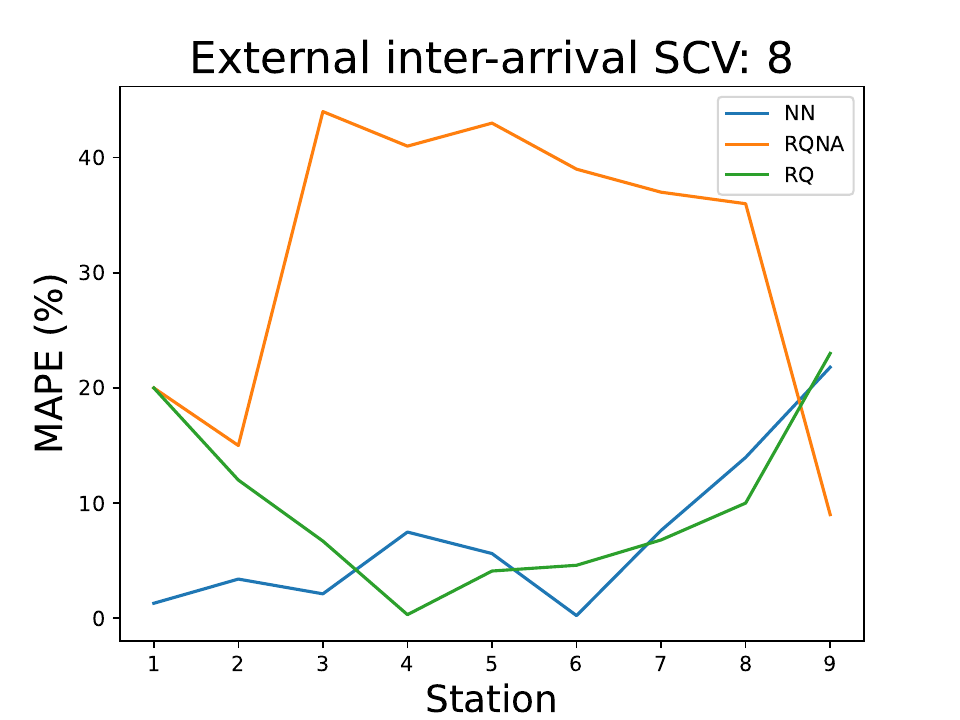}
            \caption[]%
            {{\small Hyper-exponential arrival.}}    
            \label{fig:hyper}
        \end{subfigure}
        \caption{9 Station tandem queueing system.}%
    \label{fig:mom_auto_anal}
    \end{figure}

Figures~\ref{fig:deter} and~\ref{fig:hyper} present the MAPE results for all three models with an external arrival SCV of 0 and 8, respectively. In both scenarios, our model outperforms the other models. The RQ model slightly outperforms our model in the last station with SCV 0 external arrival (see Figure~\ref{fig:deter}). The RQNA model outperforms our model in the last station, with SCV 8  external arrival (see Figure~\ref{fig:hyper}). Otherwise, our model is better in the vast majority of the cases.     

\subsection{Runtimes}\label{sec:runtimes}

In this section, we measure the runtime for each NN. The runtimes representing the inference time for each NN are tested over a personal PC, with a processor Intel(R) Core(TM) i7-14650HX   2.20 GHz, 32 GB. Table~\ref{tab:runtimes_res} presents the inference time for each NN. Each time, we examine 750 instances in parallel to demonstrate its speed advantage over simulation. As the results indicate, for 750 instances, all NNs take only a fraction of a second.  

\begin{table}[!htp]\centering
\caption{Runtime results}\label{tab:runtimes_res}
\scriptsize
\begin{tabular}{lrrr}\toprule
&Number of instances &Runtimes (Sec) \\ \hline
NN 1 &750 &0.004 \\ \hline
NN 2 &750 &0.0051 \\ \hline
NN 3 &750 &0.0045 \\ \hline
\bottomrule
\end{tabular}
\end{table}

\section{Conclusion}\label{sec:conclusions}

This paper proposes a novel machine learning-based method to analyze a tandem queueing system with an external renewal arrival process, employing three NNs. The first NN estimates the departure process of a $GI/GI/1$ queue, the second calculates the steady-state probabilities of a $G/GI/1$ queue, and the third estimates the departure process of a $G/GI/1$ queue. Our approach is unique in representing the departure process through moments and polynomial auto-correlation values.

We conducted a comprehensive performance evaluation of our NNs, and the results indicate that our method achieves an error rate of less than 5\% in most scenarios. Comparison with existing methods demonstrates that our approach is state-of-the-art.

Additionally, we provide a compelling analysis of the marginal effect of the $i^{\text{th}}$ inter-arrival and service time moments on the steady-state probabilities, showing that the first five moments predominantly determine these probabilities. Similarly, an analysis of inter-arrival auto-correlation values reveals that the first two lags of the first and second-degree polynomial auto-correlation values almost entirely determine the steady-state probabilities of a $G/GI/1$ queue. These insights could aid in developing analytical approximation methods for similar systems.



\newpage











\begin{appendices}

\section{NN architecture}\label{append:architecture}

In Table~\ref{tab:nn_archi} we present the architecture properties of \textbf{NN 1}, \textbf{NN 2} and \textbf{NN 3}. 

\begin{table}[!htp]\centering
\caption{NN architecture properties}\label{tab:nn_archi}
\scriptsize
\begin{tabular}{lrrrrrr}\toprule
&Input size &output size &\# of hidden layers &\# node of hidden layers &number of parameters \\
NN 1 &10 &18 &3 &50, 70, 50 &7925 \\
NN 2 &18 &1500 &5 &50, 70, 200, 350, 600 &1200569 \\
NN 3 &18 &18 &3 &50, 70 50 &8325 \\
\bottomrule
\end{tabular}
\end{table}

We employed the Rectified Linear Unit (ReLU) as the activation function for all hidden layers, where $ReLU(input) = max\{ input, 0\}$. Additionally, the Softmax function was applied to the output layer, defined as $SOFTMAX_i ([a_1,...,a_n]) = \frac{exp(a_i)}{\sum_{j=1}^n exp(a_j)}, i=1, \ldots, n$, ensuring that the output-layer weights sum to 1.

In terms of computational complexity, the network we used is considered small (as large networks typically contain several million parameters; for instance, see the survey in~\cite{DBLP:journals/corr/abs-1803-01164}, where NNs range from 1M to 250M parameters). This allows for increasing network complexity as our machine-learning approach is extended to more complex queueing systems. 

\section{Fine-tuning of the NN}\label{append:Fine-tuning}
We provide technical details of our hyper-params fine-tuning of our model for a given ($n_{arrival}$, $n_{service}$). The hyper-params that we wish to optimize over are learning rate, number of training epochs, number of hidden layers, number of neurons for each layer, Batch-size, and the weight-decay parameter of the Adam optimizer (this is a regularization parameter; for more details on the Adam optimizer see~\cite{8624183} and reference therein). 

The state space for each parameter is:
\begin{itemize}
    \item Learning rate: $ \{0.01, 0.05, 0.001, 0.0001\}$.
    \item Training epochs: $[100,350]$.
    \item Number of hidden layers: $[4,10]$.
    \item Number of neurons for each layer: $\{50,70,100,150,200,300\}$
    \item Batch-size: $\{64, 128, 256, 512\}$.
    \item Weight-decay: $\{10^{-4}, 10^{-5}, 10^{-6}\}$.
\end{itemize}
We do 500 random searches over these parameters and use the best settings. That is the one with the lowest SAE over the validation set.  

Each search takes an average of 0.76 hours, training on an NVIDIA GEFOrce4700 Tensor Core GPU and 8 GB memory.




\end{appendices}


\end{document}